\documentclass[preprint,12pt]{elsarticle}
\usepackage{a4wide}
\usepackage{type1cm}        % activate if the above 3 fonts are
\usepackage[dvipsnames]{xcolor}

\usepackage{makeidx}         % allows index generation
\usepackage{graphicx}        % standard LaTeX graphics tool
\usepackage{multicol} 
\usepackage{multirow}% used for the two-column index
\usepackage[bottom]{footmisc}% places footnotes at page bottom

\usepackage{bbm}                            % voor \N, \R, \C, \Q in blackboard-bold
\usepackage{csquotes}
\usepackage{newtxtext}       %
\usepackage{bm} 
\usepackage{bbm}
\usepackage{amssymb}
\usepackage{amsmath}
\usepackage{amsthm}
\theoremstyle{definition}
\usepackage{caption}
\newtheorem{thm}{Theorem}[section]
\newtheorem{lem}{Lemma}[section]

\newtheorem{rem}{Remark}[section]
\usepackage{graphicx}
\usepackage{booktabs}
\usepackage{array}
\usepackage[title]{appendix}
\usepackage{empheq}
\usepackage{subcaption}
\newcommand{\R}{\mathbbm{R}}                % reele getallen
\usepackage[normalem]{ulem}
\usepackage{algorithm}

\usepackage{hyperref}
\usepackage{algorithmic}
\makeindex             % used for the subject index
\begin{document}

\begin{frontmatter}

\title{A stabilized finite element method for a Morpho-Visco-Poroelastic Model} 

\author[label1,label2]{S. Asghar\corref{cor1}}
\ead{sabia.asghar@uhasselt.be}
\author[label3]{D.R. den Bakker}
\author[label4]{E. Javierre}
\author[label5,label6]{Q. Peng}
\author[label1,label2]{F.J. Vermolen}

\cortext[cor1]{Corresponding author}
\affiliation[label1]{organization={Computational Mathematics Group (CMAT), Department of Mathematics and Statistics, University of Hasselt},addressline={Diepenbeek}, country={Belgium}}
\affiliation[label2]{organization={Data Science Institute (DSI), University of Hasselt},addressline={Diepenbeek}, country={Belgium}}
\affiliation[label3]{organization={Delft Institute of Applied Mathematics, Delft University of Technology}, addressline={Delft}, country={The Netherlands}}
\affiliation[label4]{organization={IUMA and Applied Mathematics Department, University of Zaragoza}, addressline={Zaragoza}, country={Spain}}
\affiliation[label5]{organization={Mathematics for AI in Real-world Systems, School of Mathematical Sciences, Lancaster University}, addressline={Lancaster}, country={UK}}
\affiliation[label6]{organization={Mathematical Institute, Leiden University}, addressline={Leiden}, country={The Netherlands}}

%% Abstract
\begin{abstract}
Studying the structure of soft tissues is important and relevant in biology, particularly in some diseases, such as tumor growth and dermal contraction after burn injury. Based on the complicated characteristics of the tissue and for the sake of a better understanding of the underlying biomechanics, we propose a mathematical model that combines elastic, viscous, and porous effects with growth or shrinkage due to microstructural changes. The framework is referred to as morpho-visco-poroelasticity. Although the existence results of the solution to the problem are not given in this study, we assess the stability of the equilibria for both the continuous and semi-discrete versions of the model, and the key features of this modelling framework have been discussed. To obtain reliable numerical solutions, a stabilized finite element (FE) scheme is proposed for the morpho-visco-poroelasticity equations to avoid spurious oscillations in the pressure profile; the success of this FE scheme is verified by numerical simulations and convergence investigation in both spatial and temporal aspects. For a more quantitative assessment, the total variation of the pressure profile is evaluated as a function of the stabilization parameter. 

\end{abstract}

%% Keywords
\end{frontmatter}

\section{Introduction}
Biological tissues are complex materials consisting of a solid matrix, primarily composed of cells and the extracellular matrix (ECM), permeated by interstitial fluids that play a crucial role in their biomechanical response \cite{Simon1990}. They are frequently subjected to relatively high mechanical stresses and exhibit a non-linear relationship between stress and strain, in particular at large deformations \cite{Taylor2020}. In addition, soft tissues are capable of growth, remodeling, and adaptation \cite{Humphrey2002,Rodriguez1994}, possibly leading to permanent changes of their mechanical properties. Understanding these phenomena is essential in a variety of biomedical applications, including cartilage mechanics \cite{Mak1986}, tumor growth \cite{JainETAL2014}, wound healing \cite{Bin2011}, and tissue remodeling \cite{Humphrey2002}. In such settings, both fluid transport and solid deformation are tightly coupled and must be modeled simultaneously.

In this regard, poroelastic models have emerged as a powerful tool to characterize the mechanical behavior of fluid-saturated biological tissues. Originally developed by Biot in the context of soil mechanics \cite{Biot1941}, the theory of poroelasticity has been extensively adapted to describe biological systems such as articular cartilage \cite{Mow1980}, intervertebral discs, brain tissue, and tumor microenvironments \cite{BaxterJain1989,NettiETAL1995}. The fundamental principle of poroelasticity is that mechanical deformation and fluid transport are intrinsically and inseparably coupled. When a porous solid matrix is subjected to mechanical loading, it induces fluid flow through its interconnected pore network. In contrast, spatial and temporal variations in fluid pressure generate stresses within the solid skeleton, leading to deformation. This two-way interaction establishes a dynamic feedback mechanism between solid mechanics and fluid dynamics,
which is essential to understand the behavior of many biological tissues and plays a critical role in both physiological function and pathological conditions. As examples of where this is important, we mention the time-dependent behavior of articular cartilage under compressive loading \cite{Mow1980} or the mechanical behavior of solid tumors due to abnormal vascular permeability and deficient lymphatic drainage \cite{Wu2013}. These behaviors have a significant impact on effective drug delivery \cite{JainETAL2014,HeldinETAL2004}.
In general, poroelasticity provides a unified framework for capturing the interplay between fluid flow and mechanical deformation, allowing a deeper understanding of complex biomechanical systems.

Although poroelastic models successfully capture the coupled fluid-solid mechanics of biological tissues, they assume that the solid matrix exhibits purely elastic behavior, which may be insufficient to describe the time-dependent mechanical response of certain soft tissues. Poro-viscoelastic models extend the classical poroelastic framework by incorporating the intrinsic viscoelastic properties of the solid matrix, thus accounting for both the viscous drag associated with fluid flow through the porous medium and the inherent viscous dissipation within the solid skeleton itself \cite{Mak1986}. Mathematical models based on (visco-)poroelastic theory provide a quantitative framework to predict these coupled phenomena, offering insights that are difficult to obtain through experimental observation alone.

From a computational perspective, the numerical solution of Biot's consolidation model presents significant challenges related to the stability of finite element discretizations. Standard Galerkin formulations often exhibit spurious oscillations in the pressure field, particularly when strong pressure gradients occur or when the permeability is small relative to the mesh size \cite{Aguilar2008,BergerETAL2015}. These spurious oscillations arise from the failure of certain finite element pairs to satisfy the inf-sup (or LBB) condition uniformly with respect to the physical parameters of the problem. Several stabilization strategies have been developed to address these difficulties. One such approach involves the use of finite element spaces that satisfy the appropriate inf-sup conditions, such as the classical Taylor-Hood elements or the MINI element \cite{MuradLoula1994}. Although these methods provide theoretical stability, they may still exhibit oscillatory behavior in the presence of very sharp boundary layers or when dealing with materials of low permeability \cite{RodrigoETAL2016}. 

Aguilar et al. \cite{Aguilar2008} attribute these oscillations to a lack of monotonicity in the numerical scheme and propose a stabilization approach based on the addition of a time-dependent artificial term to the flow equation, with a stabilization parameter depending on the elastic properties of the solid and the mesh size. This stabilization parameter has been shown to be optimal in the one-dimensional case \cite{Aguilar2008}. The method provides oscillation-free solutions independently of the chosen discretization parameters while maintaining the use of linear finite elements for both displacements and pressure. 

Furthermore, Rodrigo et al. \cite{RodrigoETAL2018} introduced a stabilization approach for the three-field formulation of Biot's model based on enriching the piecewise linear continuous finite element space for displacement with edge or face bubble functions. By applying a consistent perturbation of the bilinear form, the bubble functions can be eliminated locally, resulting in a stable scheme with the same number of degrees of freedom as the classical P1-RT0-P0 discretization. This approach ensures uniform stability with respect to both physical and discretization parameters. Similarly, Berger et al. \cite{BergerETAL2015} developed a stabilized conforming mixed finite element method for the three-field formulation using the lowest possible approximation order. Their method incorporates a local pressure jump stabilization term in the mass conservation equation, leading to a symmetric linear system while ensuring stability and avoiding pressure oscillations. Moreover, divergence-conforming finite element discretizations for the Biot problem have been systematically developed and rigorously analyzed in the literature, notably by \cite{Badia2024} among others. In parallel, stress and flux reconstruction techniques have been introduced within the poroelasticity framework to enable reliable {\em a posteriori} error estimation, providing effective tools to assess and improve the accuracy of numerical approximations \cite{Riedlbeck2017}.

In addition to stability and accuracy, monotonicity plays a crucial role in the numerical treatment of coupled and transport-dominated systems. This aspect becomes particularly important in problems characterized by strong convection or chemotactic effects, where nonphysical oscillations may arise. For instance, a second-order accurate, positivity-preserving scheme has recently been proposed and analyzed for the Keller–Segel–Navier–Stokes model, demonstrating robust performance in preserving physically meaningful solutions \cite{Wang2026}.

Despite their ability to capture elastic, porous, and time-dependent responses, poroelastic and poro-viscoelastic models still neglect an important biological feature, namely microstructural evolution. Consequently, in addition to poroelastic effects, many biological tissues undergo growth or shrinkage as a result of such microstructural changes. These microstructural changes are often caused by cells that alter their direct environment. As examples, we mention the ability of cancer cells that adapt the mechanical properties in their immediate surroundings or the permanent contraction of skin following deep tissue injury. 

The morphoelasticity framework combines microstructural evolution and mechanical deformation. This model was first introduced by Rodriguez {\em et al.} \cite{Rodriguez1994} and was later further developed and systematically presented by Goriely \cite{Goriely2015}. Hall \cite{Hall2017} extended the morphoelastic formalism to multiple spatial dimensions. A limitation of the morphoelastic model is that it neglects the porous nature of biological tissues. 
To address this limitation, we extended the morphoelastic framework to address the porous and viscous nature of tissue. This is done by introducing a morpho-visco-poroelastic model. Since this is a new model, we will address some numerical properties like stability and monotonicity. The question of well-posedness in terms of existence and uniqueness is postponed to a separate study. 

The manuscript is structured as follows. The mathematical model is presented in Section \ref{sec:2}. In Section \ref{sec:3}, we discuss the stability of the model in detail. Section \ref{sec:4} presents the numerical scheme. In Section \ref{sec:5}, we present numerical results and finally we list the conclusions in Section \ref{sec:6}.

\section{The mathematical model}\label{sec:2} 
In this section, we present the system of partial differential equations that model the morpho-visco-poroelastic nature of tissues. The equations contain a momentum balance, constitutive law, continuity equation and the evolutionary equation for the effective (Eulerian) strain from \cite{Hall2017, Hall2008}.
\par The morpho-visco-poroelastic model reads as: 

\begin{subequations} \label{morphoporoelas}
  \begin{empheq}[left=\empheqlbrace]{align}
  	& \frac{D}{Dt} (\rho {\bf w}) + \rho (\nabla \cdot {\bf w}) {\bf w} - \nabla \cdot \boldsymbol{\sigma}({\bf w},\boldsymbol{\varepsilon}) + \nabla p = {\bf f}_{\bf u}, \label{mvpe_momentum}\\[1ex]
	& \boldsymbol{\sigma}({\bf w},\boldsymbol{\varepsilon}) = \boldsymbol{\sigma}_{el}(\boldsymbol{\varepsilon}) + \boldsymbol{\sigma}_{vis}({\bf w}), 
    \label{mvpe_stress} \\[1ex]
	& \nabla \cdot {\bf w} - \nabla \cdot (\kappa \nabla p) = f_p, \label{mvpe_continuity} \\[1ex]
	& \frac{D \boldsymbol{\varepsilon}}{Dt} + \boldsymbol{\varepsilon} ~ \text{skw}(\nabla {\bf w}) - \text{skw}(\nabla {\bf w}) \boldsymbol{\varepsilon}  + (\text{tr}(\boldsymbol{\varepsilon}) - 1) \text{sym}(\nabla {\bf w}) = -\mathbf{G}, \label{mvpe_strain} \\[1ex]
	& \frac{D {\bf u} }{Dt} = {\bf w}, \label{mvpe_dudt}
  \end{empheq}
\end{subequations}
in $\Omega(t)\times (0,T]$, where 
\begin{align}
    \boldsymbol{\sigma}_{el}(\boldsymbol{\varepsilon}) = 2\mu \boldsymbol{\varepsilon} + \lambda  ~\text{tr}(\boldsymbol{\varepsilon}) {\bf I}, 
    \quad
    \boldsymbol{\sigma}_{vis}({\bf w}) = \mu_1 ~ \text{sym}(\nabla {\bf w}) + \mu_2 ~ \text{tr}(\text{sym}(\nabla {\bf w})) \mathbf{I}. \tag*{(1f)}
    \label{mvpe_se_sv}
\end{align}
Here, ${\bf u}$ and $p$, respectively, represent the local displacement of the skeleton and the pore pressure, ${\bf w}$ the displacement velocity, $\boldsymbol{\varepsilon}$ the strain tensor, ${\bf G}$ the net growth tensor, and $\kappa$ represents the ratio of the permeability of the porous medium and the fluid viscosity. Further, $D(\cdot)/Dt$ denotes the material derivative, and $\text{sym}(.)$ and $\text{skw}(.)$, respectively, denote the symmetric and skew-symmetric part of a matrix or tensor.  Eq. (\ref{mvpe_momentum}) represents the balance of momentum where inertial terms are taken into account, $\rho$ denotes the tissue density, $\boldsymbol{\sigma}$ denotes the effective stress tensor and the right-hand term ${\bf f}_{\bf u}$ represents applied body force(s). Eq. (\ref{mvpe_stress}) accounts for the elastic and viscoelastic behavior of the material, with their respective constitutive relations given in \ref{mvpe_se_sv}. There $\lambda$ and $\mu$ denote the Lam\'e parameters, $\mu_1$ and $\mu_2$ the viscosity coefficients and $\mathbf{I}$ the identity tensor. Eq. (\ref{mvpe_continuity}) represents the conservation equation for the fluid, where Darcy's law \cite{Hubbert1956},
\begin{align*}
{\bf v} + \kappa ~ \nabla p = 0,
\end{align*}
has already been incorporated. Here, ${\bf v}$ represents the fluid flow velocity. In Eq. (\ref{mvpe_continuity}) the source term $f_p$ represents a forced fluid extraction or injection process. Eq. (\ref{mvpe_strain}) represents the evolution equation of the strain tensor. Note that the primary variables in Eqs. (\ref{morphoporoelas}) are ${\bf w}$, $\boldsymbol{\varepsilon}$ and $p$. The displacement ${\bf u}$ can be obtained by integrating ${\bf w}$ over one time-step as a post-processing step.

Note that Eq. (\ref{mvpe_strain}) models the evolution of tissue strain as a result of microstructural changes. This will lead to permanent (local) displacements, as the ${\bf G}$-tensor quantifies growth or shrinkage of the medium as a result of microstructural changes. If ${\bf G} = {\bf 0}$, then the model does not predict any permanent displacements. 

The model needs to be completed with appropriate boundary and initial conditions. Unless stated otherwise, we consider the following:
\begin{subequations} \label{mpve_ibc}
  \begin{empheq}[left=\empheqlbrace]{align}
   & {\bf w} = \mathbf{0}, \quad \kappa \nabla p \cdot \boldsymbol{n} = g_N , & \text{on $\Gamma_1(t) \times (0,T]$,} \label{mpve_bc1} \\[1ex]
   & \boldsymbol{\sigma}\cdot\boldsymbol{n} = \mathbf{f}_b, \quad p = p_0, & \text{on $\Gamma_2(t) \times (0,T]$,}  \label{mpve_bc2} \\[1ex]
   & {\bf w}({\bf x},0) = {\bf 0}, & \text{in $\Omega(0)$,} \label{mpve_ic1} \\[1ex]
   & {\boldsymbol \varepsilon}({\bf x},0) = {\boldsymbol 0}, & \text{in $\Omega(0)$,} \label{mpve_icw} 
  \end{empheq}
\end{subequations}
where $g_N$ represents the fluid of the flow over $\Gamma_1$, ${\bf f}_b$ represents an external force on the boundary, and $p_0$ denotes the environmental atmospheric pressure, which is taken to be constant.
We note that if $\kappa$ is a constant (not depending on the position or time) and if the fluid flow velocity stays finite, then from Darcy's Law heuristically it follows that $\nabla p \longrightarrow {\bf 0}$ as $\kappa \longrightarrow \infty$. This implies that $p = c$, for some $c \in \R$. Since $p = p_0$ on $\Gamma_2(t)$, it follows that $p = p_0$ in $\Omega$.
It should be noted that the standard morphoelastic formulation \cite{EgbertsStab2021} is recovered when $\kappa \to +\infty$ in Eqs. \eqref{mpve_ibc}. Furthermore, it is important to mention that the proposed model is formulated in dimensional form, where all variables and parameters retain their physical units (the nomenclature has been provided in Appendix \ref{Appendix}).

\section{Stability of the model} \label{sec:3}
The term "stable" informally means resistant to change. For technical use, the term has to be defined
more precisely in terms of the mathematical model, but the same connotation applies. Stability analysis is crucial to ensuring the reliability and accuracy of the numerical methods used to solve these problems. In this section, we discuss the stability of the morpho-visco-poroelastic model. We first prove the symmetry of the strain tensor which further leads to the stability of the model. Then we assess the \textit{linear} stability of the steady state problems using the Fourier series.
\subsection{Symmetry of the strain tensor}
We demonstrate that the strain tensor will remain symmetric at all times if it is initially symmetric. Before we prove the result on the symmetry of $\boldsymbol{\varepsilon}$, recall the following operations with tensors. If $A,B\in \mathbb{R}^{n\times n}$, then 

\begin{align}
    A : B = \sum_{i,j = 1}^n A_{ij}B_{ij}.
    \label{ten-scal-prod1}
\end{align}
Multiplication of matrices $A^T$ and $B$ gives componentwise
$$\displaystyle (A^T B)_{ij} = \sum_{k=1}^n A_{ki} B_{kj},$$ and consequently
\begin{align}
	\text{tr}(A^T B) = \sum_{i=1}^n  (A^T B)_{ii} =  \sum_{i=1}^n \sum_{k=1}^n A_{ki} B_{ki} = A : B.
\label{ten-scal-prod2}
\end{align}

\begin{lem}\label{lemma1}
Let ${\boldsymbol v}, {\boldsymbol L} \in \mathbb{R}^{d \times d}$ be d-dimensional tensors. Suppose that ${\boldsymbol L}$ is skew-symmetric (${\boldsymbol L}^T = -{\boldsymbol L}$), then for all ${\boldsymbol v} \in \mathbb{R}^{d \times d}$ the tensorial scalar product satisfies
\begin{align*}
{\boldsymbol v} : ({\boldsymbol L} {\boldsymbol v}) = {\boldsymbol v} : ({\boldsymbol v} {\boldsymbol L} ) = 0. 
\end{align*}
\end{lem}
\noindent {\bf Proof.} 
Choose any $\boldsymbol{v} \in \mathbb{R}^{d \times d}$, use Eq. (\ref{ten-scal-prod2}) and ${\boldsymbol L}^T = -{\boldsymbol L}$: %, $(AB)^T = B^T A^T$, and equation (\ref{ten-scal-prod2}), then we arrive at
\begin{align*}
\boldsymbol{v} : (\boldsymbol{L} \boldsymbol{v}) 
=\text{tr}(\boldsymbol{v}^T \boldsymbol{L}\boldsymbol{v})
=-\text{tr}(\boldsymbol{v}^T\boldsymbol{L}^T\boldsymbol{v})
=-\text{tr}((\boldsymbol{L}\boldsymbol{v})^T\boldsymbol{v})
=-\boldsymbol{v} : (\boldsymbol{L}\boldsymbol{v}).
\end{align*}
Hence $\displaystyle{ {\boldsymbol{v}} : ({\boldsymbol L} {\boldsymbol{v}}) = 0 }.$ \\

Furthermore, since $A : B = A^T : B^T$, take $\boldsymbol{w} = \boldsymbol{v}^T$ and we use ${\boldsymbol L}^T=-{\boldsymbol L}$: %, we get using the above relation:
\begin{align*}
(\boldsymbol{w} \boldsymbol{L}) : \boldsymbol{w} =  (\boldsymbol{w} \boldsymbol{L})^T : \boldsymbol{w}^T = (\boldsymbol{L}^T \boldsymbol{v}) : \boldsymbol{v} = -(\boldsymbol{L} \boldsymbol{v}) : \boldsymbol{v} = -\boldsymbol{v} : (\boldsymbol{L} \boldsymbol{v})  = 0.
\end{align*}
This proves the lemma. \hfill $\Box$ \\[2ex]
Note that the above products can be generalized to scalar products of anti-Hermitian operators on Hilbert spaces. Using this fact, the only thing that remained to be done in an alternative proof was to demonstrate that ${\boldsymbol v} : ({\boldsymbol w} \boldsymbol L)$ was a proper inner product. Based on the above lemma, we demonstrate the following claim:
\begin{thm}\label{Th_aux}
Let $\bf{G}$ be a symmetric tensor for $t\geq 0$ and $\boldsymbol{\varepsilon}$ satisfy
\begin{align}
    \frac{D \boldsymbol{\varepsilon}}{Dt} + \boldsymbol{\varepsilon} ~ \text{skw}(\nabla {\bf w}) - \text{skw}(\nabla {\bf w}) \boldsymbol{\varepsilon}  + (\text{tr}(\boldsymbol{\varepsilon}) - 1) \text{sym}(\nabla {\bf w}) = -{\bf{G}}, 
    \label{eps}
\end{align}
in an open Lipschitz domain $\Omega\in\mathbb{R}^d$ for $t>0$. Suppose that $\boldsymbol{\varepsilon}$ is symmetric on $t=0$, then $\boldsymbol{\varepsilon}$ remains symmetric for $t>0$. 
\end{thm}
\noindent {\bf Proof}. We take the transpose of Eq. (\ref{eps}) and use $$\text{skw}(\nabla{\bf w})^T=-\text{skw}(\nabla{\bf w}) ~\text{ and }~ \text{sym}(\nabla{\bf w})^T=\text{sym}(\nabla{\bf w}).$$ Then we obtain the following
\begin{align}
	\frac{D \boldsymbol{\varepsilon}^T}{Dt} + \boldsymbol{\varepsilon}^T \text{skw}(\nabla {\bf w})  - \text{skw}(\nabla {\bf w}) \boldsymbol{\varepsilon}^T + (\text{tr}(\boldsymbol{\varepsilon}) - 1) \text{sym}(\nabla {\bf w}) = -\bf{G}^T.
	\label{eps_T}
\end{align}
Subtraction of Eqs. (\ref{eps}) and (\ref{eps_T}) gives
\begin{align*}
	\frac{D \left( \boldsymbol{\varepsilon} - \boldsymbol{\varepsilon}^T \right)}{Dt} + \left( \boldsymbol{\varepsilon} - \boldsymbol{\varepsilon}^T \right)  ~ \text{skw}(\nabla {\bf w}) - \text{skw}(\nabla {\bf w}) \left( \boldsymbol{\varepsilon} - \boldsymbol{\varepsilon}^T \right)  = - \left( \bf{G} - \bf{G}^T \right).
\end{align*}
If the growth tensor $\bf G$ remains symmetric on time, then we have
\begin{align}
	\frac{D \left( \boldsymbol{\varepsilon} - \boldsymbol{\varepsilon}^T \right)}{Dt} + \left( \boldsymbol{\varepsilon} - \boldsymbol{\varepsilon}^T \right)  ~ \text{skw}(\nabla {\bf w}) - \text{skw}(\nabla {\bf w}) \left( \boldsymbol{\varepsilon} - \boldsymbol{\varepsilon}^T \right)  = 0.
	\label{aux1}
\end{align}
From this equation, it is clear that $\boldsymbol{\varepsilon} - \boldsymbol{\varepsilon}^T = 0$ represents an equilibrium solution. In addition, we can prove that $\boldsymbol{\varepsilon} - \boldsymbol{\varepsilon}^T = 0$ is the only solution provided that $\boldsymbol{\varepsilon} - \boldsymbol{\varepsilon}^T = 0$ is fulfilled at $t=0$. \\

\noindent Taking $\boldsymbol{v} = \boldsymbol{\varepsilon} - \boldsymbol{\varepsilon}^T $ and $\boldsymbol{L} = \text{skw}(\nabla {\bf w})$, Eq. (\ref{aux1}) can be rewritten as
\begin{align*}
	\frac{D \boldsymbol{v} }{Dt} + \boldsymbol{v}  ~ \boldsymbol{L} - \boldsymbol{L} ~ \boldsymbol{v}  = 0.
\end{align*}
Then we have that 
\begin{align*}
	\boldsymbol{v} : \frac{D \boldsymbol{v} }{Dt} + \boldsymbol{v} : \left( \boldsymbol{v} \boldsymbol{L} \right) - \boldsymbol{v} : \left( \boldsymbol{L} \boldsymbol{v} \right) = 0.
\end{align*}
Using Lemma \ref{lemma1}, the previous equation reduces to 
\begin{align*}
	\boldsymbol{v} : \frac{D \boldsymbol{v} }{Dt} = 0.
\end{align*}
Define $||\boldsymbol{v}||^2 := \boldsymbol{v} : \boldsymbol{v}$, then it follows that 
\begin{equation*}
\frac12 \frac{{D}}{{D}t} \, ||\boldsymbol{v}||^2 = 0,
\end{equation*}
from which we obtain that $||\boldsymbol{v}({\bf x},t)||^2= ||\boldsymbol{v}({\bf x},0)||^2$ for any $t>0$. This proves that if 
$\boldsymbol{\varepsilon}$ is symmetric on $t=0$ (i.e. $\boldsymbol{v}({\bf x},0) = 0$) then $\boldsymbol{\varepsilon}$ is symmetric for $t>0$. 
\hfill $\Box$ 

The previous result can be further generalized to a broader class of growth tensors. 
Further, we demonstrate that small perturbations around symmetric strain tensor remain small, which is a characteristic of stability.
\begin{thm}\label{Theorem5}
Let $\boldsymbol{\varepsilon}$ satisfy
\begin{align*}
    \frac{D \boldsymbol{\varepsilon}}{Dt} + \boldsymbol{\varepsilon} ~ \text{skw}(\nabla {\bf w}) - \text{skw}(\nabla {\bf w}) \boldsymbol{\varepsilon}  + (\text{tr}(\boldsymbol{\varepsilon}) - 1) \text{sym}(\nabla {\bf w}) = -\alpha \boldsymbol{\varepsilon}, 
\end{align*}
in an open Lipschitz domain $\Omega\in\mathbb{R}^d$ for $t>0$. Suppose that $\boldsymbol{\varepsilon}$ is symmetric for $t=0$, then $\varepsilon$ remains symmetric for $t>0$, and symmetry is asymptotically stable with respect to perturbations if and only if $\alpha>0$ (and stable if and only if $\alpha\geq 0$).
\end{thm}

\noindent {\bf Proof}. Reproducing the steps of the proof of Theorem \ref{Th_aux} we obtain the following equation for $\boldsymbol{v} = \boldsymbol{\varepsilon} - \boldsymbol{\varepsilon}^T $:
\begin{align*}
	\frac{D \boldsymbol{v} }{Dt} + \boldsymbol{v}  ~ \boldsymbol{L} - \boldsymbol{L} ~ \boldsymbol{v}  = -\alpha \boldsymbol{v},
\end{align*}
where $\boldsymbol{L} = \text{skw}(\nabla {\bf w})$. Using the tensor product and Lemma \ref{lemma1} we obtain that 
\begin{equation*}
\boldsymbol{v} : \frac{{D}}{{D}t} \, \boldsymbol{v} = -\alpha\,\boldsymbol{v} : \boldsymbol{v},
\end{equation*}
and 
\begin{equation*}
\frac12 \frac{{D}}{{D}t} \, ||\boldsymbol{v}||^2 = -\alpha \, ||\boldsymbol{v}||^2.
\end{equation*}
Integrating over time from $t = 0$ and using $\boldsymbol{v} = 0$ at $t = 0$ and $\alpha \geq 0$, gives 
\begin{equation}
0 \le \frac{1}{2} ||\boldsymbol{v}||^2 = - \int_0^t \alpha \, ||\boldsymbol{v}||^2 \, \mathrm{d}s \le 0.
\end{equation}

This implies that $||\boldsymbol{v}|| = 0$ on $t > 0$ if $||\boldsymbol{v}|| = 0$ on $t = 0$. Hence $\boldsymbol{v} = 0$ for $t > 0$, which represents symmetry, is the only possibility if $\boldsymbol{v} = 0$ on $t = 0$. Furthermore, stability is also proved by this argument.
\hfill $\Box$
\subsection{Stability of the steady state}
Consider the one-dimensional version of Eqs. \eqref{morphoporoelas} in the fixed domain $\Omega = (0,1)$. If we further assume $\rho$ to remain constant, $G = \alpha \varepsilon$ ($\alpha \ge 0$), $f_u = 0$ and $f_p = 0$, we obtain:

\begin{subequations} \label{mpve_1D}
  \begin{empheq}[left=\empheqlbrace]{align}
  	& \rho \left [\frac{\partial w}{\partial t} + 2w \frac{\partial w}{\partial x} \right ]- E \frac{\partial \varepsilon}{\partial x} - \mu_{vis} \frac{\partial^2 w}{\partial x^2} + \frac{\partial p}{\partial x}= 0, \\%[1ex]
	& \frac{\partial w}{\partial x} - \frac{\partial}{\partial x}\left(\kappa \frac{\partial p}{\partial x}\right) = 0, \\[0.5ex]
	& \frac{\partial \varepsilon}{\partial t} + (\varepsilon - 1) \frac{\partial w}{\partial x} = -\alpha \varepsilon, \\%[1ex]
	& \frac{\partial u }{\partial t} = w, 
  \end{empheq}
\end{subequations}
 
where $E = 2\mu+\lambda$ and $\mu_{vis} = \mu_1+\mu_2$. 
For simplicity, we take the following boundary and initial conditions:
\begin{subequations} \label{mpve_ibc_1D}
  \begin{empheq}[left=\empheqlbrace]{align}
   & w = 0, \quad \kappa \frac{\partial p}{\partial x} = 0, & \text{on $\Gamma_1 \times (0,T]$,} \\%[1ex]
   & \frac{\partial \sigma}{\partial x} = f_b, \quad p = p_0, & \text{on $\Gamma_2 \times (0,T]$,}   \\%[1ex]
   & w(x,0) = 0, & \text{in $(0,1)$,}  \\%[1ex]
   & \varepsilon(x,0) = 0, & \text{in $(0,1)$.} %\label{mpve_icw} 
  \end{empheq}
\end{subequations}
The steady state solutions of the one-dimensional model, Eqs. (\ref{mpve_1D})-(\ref{mpve_ibc_1D}) are $(\overline{w},\overline{\varepsilon},\overline{p}) = (0,0,p_0)$, where $p_0$ was introduced earlier as a constant. 

In the following subsections, we analyze the linear stability of the steady states of the one-dimensional model. The analysis that we perform is commonly known as Von Neumann stability analysis in cases of equidistant meshes and constant coefficients, which is based on decomposition of motion into normal modes, often using Fourier analysis, and superposition. The Von Neumann stability analysis provides sufficient and necessary conditions for linear numerical stability \cite{Fletcher1998}. The analysis looks at the growth or
decay of perturbations between consecutive time-steps, and can be implemented using standard linear algebraic
procedures. A more severe restriction is that it strictly applies only to linear systems. Despite this limitation, it is frequently applied to nonlinear systems through linearization.

\subsubsection{Stability of steady states in the continuous problem}
In this subsection, we analyze the {\em linear} stability of the one-dimensional model, Eqs. (\ref{mpve_1D})-(\ref{mpve_ibc_1D}) %(\ref{morphoporoelas}) 
with respect to sinusoidal perturbations around the equilibrium using Fourier series. We do this analysis in order to understand the \textit{a-priori} behavior of the solution. 
\par {From now on, we make the following considerations:
\begin{itemize}
\item For simplicity, we take $\Omega = (0,1)$. 
\item The stability conditions will be formulated in terms of the input parameters. 
\item We indicate the equilibria by overlines and the perturbations by hats and proceed as in \cite{EgbertsStab2021}.
\end{itemize}}

\noindent Let us consider the following perturbations %(denoted by hats) 
of the steady state solutions
\begin{equation}
{w=\overline{w}+\hat{w}=\hat{w}},\quad \varepsilon=\overline{\varepsilon}+\hat{\varepsilon}=\hat{\varepsilon},\quad p=\overline{p}+\hat{p} = p_0+\hat{p}.
\end{equation}

The linearized equations around the equilibria ${\overline{w}}=0,\ \overline{\varepsilon}=0$, and $\overline{p} = p_0$, are as follows

\begin{equation}
\left\{
\begin{array}{ll}
\displaystyle{ \rho \frac{\partial\hat{w}}{\partial t}-\mu_{vis}\frac{\partial^2 \hat{w}}{\partial x^2}-E\frac{\partial\hat{\varepsilon}}{\partial x}+\frac{\partial\hat{p}}{\partial x}=0,} \\[1.5ex]
\displaystyle{\frac{\partial\hat{w}}{\partial x} =\kappa \frac{\partial^2 \hat{p}}{\partial x^2},} \\[1.5ex]
\displaystyle{\frac{\partial\hat{\varepsilon}}{\partial t}-\frac{\partial\hat{w}}{\partial x}=-\alpha \hat\varepsilon,} \\[1.5ex]
\end{array}
\right. 
\label{S1Dmorphoporoelas}
\end{equation}
where %the variations 
$\hat{w}, \hat{\varepsilon}\ \text{and}\ \hat{p}$ are written in terms of complex Fourier series:
\begin{equation}
 f(x,t)= \sum_{k=-\infty}^{\infty}{\hat{f}_{k}(t) e^{2 i \pi k x}},
 \label{Var}
\end{equation}
with $i$ being the imaginary unit.

For the computations in the proof of this analysis, see \cite{Asghar12025}. Linear stability is guaranteed if and only if 
\begin{align}
    4\pi^2 l^2  (\alpha \mu_{vis} + E) + \dfrac{\alpha}{\kappa} \ge 0, \qquad 4\pi^2 l^2 \mu_{vis} + \dfrac{1}{\kappa} + \alpha \rho \ge 0,
\end{align}
which implies linear stability for all modes if and only if $\alpha \ge 0$.
The findings of the analysis are summarized below:
\\[1ex]

\begin{thm}\label{Theorem6}
Let $\alpha \ge 0$, then the equilibria $(\overline w,\overline{\varepsilon},\overline p) = (0,0,p_0)$ in the one-dimensional (continuous) morpho-visco-poroelastic model are linearly stable.
\end{thm}

\subsubsection{Stability of steady states in the semi-discrete problem}

Next, we analyze the stability of the semi-discrete problem under the same consideration as above. Let $h$ be the gridsize, taken uniform, and $k$ denote the index of the meshpoint positioned at $k h$. We consider $n$ nodes in total, hence $n h = 1$. The finite difference method (FDM), based on central differences for the first- and second-order spatial derivatives, gives:

\begin{equation}
\left\{
\begin{array}{ll}
\displaystyle{{w_k}^{'}(t) = \frac{\mu_{vis}}{\rho} \frac{w_{k-1}-2 w_k+w_{k+1}}{h^2}+\frac{E}{\rho}\frac{{\varepsilon_{k+1}-\varepsilon_{k-1}}}{2 h}-\frac{1}{\rho}\frac{p_{k+1}-p_{k-1}}{2 h},} \\[1.5ex]
\displaystyle{\frac{w_{k+1}-w_{k-1}}{2 h}-\kappa\frac{p_{k-1}-2 p_k+p_{k+1}}{h^2}=0,} \\[1.5ex]
\displaystyle{{\varepsilon_k}^{'}(t) = \frac{w_{k+1}-w_{k-1}}{2 h}-\alpha\varepsilon_k.} \\[1.5ex]
\end{array}
\right. 
\label{Dis1D}
\end{equation}
Let 
\begin{equation}
 w_k= \sum_{j=1}^{n-1}{\hat{w}_{j}(t) e^{-2 i \pi k j h}}, 
 \quad
 \varepsilon_k= \sum_{j=1}^{n-1}{\hat{\varepsilon}_{j}(t) e^{-2 i \pi k j h}}, 
 \quad 
 p_k= \sum_{j=1}^{n-1}{\hat{p}_{j}(t) e^{-2 i \pi k j h}}.
 \label{VarDist1}
\end{equation}
Substitution of the variations (\ref{VarDist1}) into the simplified linearised equations (\ref{Dis1D}) gives

\begin{equation}
\left\{
\begin{array}{ll}
\displaystyle{\sum_{j=1}^{n-1} {\hat{w}_j^{'}(t)}= \frac{\mu_{vis}}{h^2\rho}\sum_{j=1}^{n-1} \hat{w}_j(t) \{e^{-2 i \pi (k-1) j h }-2 e^{-2i\pi kjh}+ e^{-2i\pi (k+1)jh}\}}\\[1.5ex]
+\displaystyle{\frac{E}{2h\rho}\sum_{j=1}^{n-1} \hat{\varepsilon}_j(t) \{e^{-2i \pi (k+1)jh }-{e^{-2i\pi (k-1)jh}}\}
-\frac{1}{2h\rho}\sum_{j=1}^{n-1} \hat{p}_j(t) \{e^{-2i \pi (k+1)jh }-{e^{-2i\pi (k-1)jh}}\}},\\[1.5ex]
\displaystyle{\frac{1}{2h}\sum_{j=1}^{n-1} \hat{w}_j(t) \{e^{-2 i \pi (k+1) j h }- e^{-2i\pi (k-1)jh}\}-\frac{\kappa}{h^2}\sum_{j=1}^{n-1} \hat{p}_j(t) \{e^{-2 i \pi (k-1) j h }-2 e^{-2i\pi kjh}+ e^{-2i\pi (k+1)jh}\}=0}, \\[1.5ex]
\displaystyle{\sum_{j=1}^{n-1}{\hat{\varepsilon}_j}^{'}(t) = \frac{1}{2h}\sum_{j=1}^{n-1} \hat{w}_j(t) \{e^{-2 i \pi (k+1) j h }-e^{-2i\pi (k-1)jh}\}-\alpha \sum_{j=1}^{n-1} \hat{\varepsilon}_j(t) e^{-2 i \pi k j h }.} \\[1.5ex]
\end{array}
\right. 
\label{Dis1DFourier}
\end{equation}
Multiplication of the above set of equations by $ e^{2 i \pi k l h }$, $l\in \{1,\ldots,n-1\}$, using Euler's Formula and $2\cos{(2 \pi l h)-2}=-4 \sin^2{(\pi l h)}$ results in

\begin{equation}
\left\{
\begin{array}{ll}
\displaystyle{ {\hat{w}_l}^{'}(t) = -\frac{4\mu_{vis}}{h^2\rho} \sin^2{(\pi l h)}} ~ \hat{w}_l(t)
-\frac{i E}{h\rho}  \sin{(2 \pi l h)} ~ \hat{\varepsilon}_l(t)
+\frac{i}{h\rho}  \sin{(2 \pi l h)} ~ \hat{p}_l(t),\\[1.5ex]
\displaystyle{-\frac{i}{h}  \sin{(2 \pi lh)}~\hat{w}_l(t)+\frac{4 \kappa}{h^2} \sin^2{(\pi l h)}~\hat{p}_l(t)=0}, \\[1.5ex]
\displaystyle{{\hat{\varepsilon}_l}^{'}(t) =-\frac{i}{h} \sin{(2 \pi l h)}~\hat{w}_l(t)-\alpha ~ \hat{\varepsilon}_l(t).} \\[1.5ex]
\end{array}
\right. 
\label{Dis1Dsys}
\end{equation}
After incorporating the value of $\hat{p}_{l}(t)=\frac{ih}{2\kappa} \cot{(\pi l h)} ~ \hat{w}_{l}$ in the above system, we get

\begin{equation}
\setlength{\arraycolsep}{2em}
\frac{d}{dt}
\begin{pmatrix}
    \hat{w}_{l}(t) \\ \\[1ex]  \hat{\varepsilon}_{l}(t)\\
\end{pmatrix} +  
   \begin{pmatrix}
    \frac{4\mu_{vis}}{h^2 \rho}\sin^2{(\pi lh)}+\frac{1}{\kappa \rho}\cos^2{(\pi lh)} &\frac{iE}{h \rho}{\sin{(2 \pi lh)}} \\ \\[1ex]
    \frac{i}{h}\sin{(2\pi lh)} & \alpha \\
\end{pmatrix} 
\begin{pmatrix}
    \hat{w}_{l}(t) \\ \\[1ex] \hat{\varepsilon}_{l}(t)\\
\end{pmatrix} = \bf{0}.\\[1.5ex]
\end{equation}

The fact that the eigenvalues of the above matrix need to be non-negative (or have a non-negative real part), leads to the following linear stability criteria:
\begin{align} 
    \!\!\! \frac{4\sin^2 (\pi l h)}{h^2} \left[ \mu_{vis} \alpha + E \cos^2(\pi l h)\right] + \frac{\alpha}{\kappa} \cos^2(\pi l h) \ge 0,
    \quad
    \frac{4\mu_{vis}}{h^2} \sin^2(\pi l h) + \frac{1}{\kappa} \cos^2(\pi l h) + \alpha \rho \ge 0.
    \label{eq:stab_discrete}
\end{align}
Hence, we conclude:
\begin{thm}\label{Theorem7}
{
Let $\alpha \ge 0$, then the equilibria $(\overline w,\overline \varepsilon,\overline p) = (0,0,p_o)$ in the one-dimensional semi-discrete morpho-visco-poroelastic model are linearly stable. 
}
\end{thm}
\begin{rem}
Recall that 
\begin{equation*}
\lim_{h \rightarrow 0} \frac{4\sin^2{(\pi lh)}}{h^2}=4 \pi^2 l^2.
\end{equation*}
Hence, Theorem \ref{Theorem7} is consistent with the stability analysis of the continuous model, Theorem \ref{Theorem6}, for $h\to 0$.
\end{rem}
\begin{rem}
Applying, for instance, an implicit Euler time integration method, will always lead to a stable numerical solution.
\end{rem}

\section{The numerical scheme} \label{sec:4}
For the sake of presentation, we are dealing with a problem in two spatial dimensions, $d=2$. Then, we can write 
\[
    {\bf w} = \begin{pmatrix} w^1 \\ w^2 \end{pmatrix}, 
    \qquad
    \boldsymbol{\varepsilon} = \begin{pmatrix} \varepsilon^{11} & \varepsilon^{12} \\ \varepsilon^{21} & \varepsilon^{22} \end{pmatrix}.
\]
Taking into the account the initial conditions and Theorem \ref{Theorem5}, the tensor $\boldsymbol{\varepsilon}$ is symmetric. Furthermore, we limit the presentation to the case ${\bf G}=\alpha \boldsymbol{\varepsilon}$ in Eq. (\ref{mvpe_strain}).

\subsection{Galerkin approximation}
In this section we consider the finite element discretization of the morpho-visco-poroelastic model, Eqs. (\ref{morphoporoelas})-(\ref{mpve_ibc}). 
In order to obtain the weak formulation of the problem, for each $t\in (0,T]$ we consider the following functional spaces:
\begin{align*}
	& \mathcal{W}_t = \left\{ \mathbf{w}(\cdot,t) \in \Big( H^1(\Omega(t)) \Big)^d \,\, \Big| \,\, \mathbf{w}(\cdot, t) = \mathbf{0} \text{  on  } \Gamma_1(t) \right\}, \\
	& \mathcal{Q}_t = \left\{ p(\cdot,t) \in H^1(\Omega(t)) \,\, \mid \,\, p = p_0 \text{  on  } \Gamma_2(t) \right\}, \\ 
	& \mathcal{Q}_{0,t} = \left\{ p(\cdot,t) \in H^1(\Omega(t)) \,\, \mid \,\, p = 0 \text{  on  } \Gamma_2(t) \right\}, \\
	& \mathcal{E}_t = \left\{ \boldsymbol{\varepsilon}(\cdot,t) \in \Big( L^2(\Omega(t)) \Big)^{d\times d} \right\}.    
\end{align*}

Applying Reynold's Transport Theorem \cite{Reynolds1903}, the weak formulation of Eqs. (\ref{morphoporoelas})-(\ref{mpve_ibc}) reads: for each $t\in (0,T]$, find $(\mathbf{w}(t), p(t), \boldsymbol{\varepsilon}(t)) \in \mathcal{W}_t \times \mathcal{Q}_t \times \mathcal{E}_t$ such that 

\begin{subequations} \label{mvpe_weakform}
\begin{align}
    & \rho \frac{d}{dt} \int_{\Omega(t)} {\bf w} \cdot {\boldsymbol \phi} d \Omega + \int_{\Omega(t)} \left( \boldsymbol{\sigma}_{el}(\boldsymbol{\varepsilon}) + \boldsymbol{\sigma}_{vis}(\bf w) \right) : \text{sym}(\nabla {\boldsymbol \phi}) d \Omega - \int_{\Omega(t)} p \nabla \cdot {\boldsymbol \phi} d \Omega \nonumber \\
	& \qquad\qquad = \int_{\Omega(t)} {\bf f}_{\bf u} \cdot {\boldsymbol \phi} d \Omega + \int_{\Gamma_2(t)} {\boldsymbol f}_b \cdot {\boldsymbol \phi} d \Gamma - \int_{\Gamma_2(t)} p_0 \boldsymbol{\phi}\cdot \boldsymbol{n} d\Gamma, 
    \label{wf_w}
    \\[1ex]
    & \frac{d}{dt}\int_{\Omega(t)} \boldsymbol{\varepsilon}:\boldsymbol{\zeta} d\Omega 
    +
    \int_{\Omega(t)} \left[ \boldsymbol{\varepsilon}\,\text{skw}(\nabla {\bf w}) - \text{skw}(\nabla {\bf w})\,\boldsymbol{\varepsilon} 
    + 
    \left( \text{tr}(\boldsymbol{\varepsilon})-1 \right)\,\text{sym}(\nabla {\bf w})
    - 
    (\nabla \cdot {\bf w}) \boldsymbol{{\varepsilon}}
    \right] : \boldsymbol{\zeta} ~ d\Omega \nonumber \\
    & \qquad\qquad + \int_{\Omega(t)} \alpha \boldsymbol{\varepsilon} : \boldsymbol{\zeta} d\Omega = 0 
    \label{wf_varepsilon}
    \\[1ex]
    & \int_{\Omega(t)} \psi \nabla \cdot {\bf w} d \Omega + \int_{\Omega(t)} \kappa \nabla p \cdot \nabla \psi d \Omega = \int_{\Omega(t)} f_p \psi d\Omega - \int_{\Gamma_1(t)} g_N \psi d\Gamma 
    \label{wf_p}
\end{align}
\end{subequations}
for all $(\boldsymbol{\phi}, \psi, \boldsymbol{\zeta}) \in \mathcal{W}_t \times \mathcal{Q}_{0,t} \times \mathcal{E}_t$.
Prior to obtain the Galerkin finite element approximation of the weak formulation, Eqs. \eqref{mvpe_weakform}, let us work out the central term in Eq. \eqref{wf_varepsilon}. 
\begin{rem} \label{rem:wf_eps}
Let us consider the following choices of $\boldsymbol{\zeta} \in \mathcal{E}_t$: 
\begin{itemize}
    \item If $\boldsymbol{\zeta} = \begin{pmatrix} \zeta & 0 \\ 0 & 0 \end{pmatrix}$, then 
    \begin{align*}
        &  \int_{\Omega(t)} \left[ \boldsymbol{\varepsilon}\,\text{skw}(\nabla {\bf w}) - \text{skw}(\nabla {\bf w})\,\boldsymbol{\varepsilon} 
        + 
        \left( \text{tr}(\boldsymbol{\varepsilon})-1 \right)\,\text{sym}(\nabla {\bf w})
        - 
        (\nabla \cdot {\bf w}) \boldsymbol{{\varepsilon}} \right] : \boldsymbol{\zeta} ~ d\Omega = \\
        & \qquad
        = \int_{\Omega(t)} \left[ -\varepsilon^{11} \dfrac{\partial w^2}{\partial y} + \varepsilon^{12} \left( \frac{\partial w^2}{\partial x} - \frac{\partial w^1}{\partial y} \right) + \varepsilon^{22} \frac{\partial w^1}{\partial x} \right] \zeta ~d\Omega
        -
        \int_{\Omega(t)} \frac{\partial w^1}{\partial x} \zeta ~d\Omega.
    \end{align*}
    \item If $\boldsymbol{\zeta} = \begin{pmatrix} 0 & \zeta \\ 0 & 0 \end{pmatrix}$, then 
   
    \begin{align*}
        &  \int_{\Omega(t)} \left[ \boldsymbol{\varepsilon}\,\text{skw}(\nabla {\bf w}) - \text{skw}(\nabla {\bf w})\,\boldsymbol{\varepsilon} 
        + 
        \left( \text{tr}(\boldsymbol{\varepsilon})-1 \right)\,\text{sym}(\nabla {\bf w})
        - 
        (\nabla \cdot {\bf w}) \boldsymbol{{\varepsilon}} \right] : \boldsymbol{\zeta} ~ d\Omega = \\
        & \qquad
        = \int_{\Omega(t)} \left[ \varepsilon^{11}\frac{\partial w^1}{\partial y} - \varepsilon^{12}\left( \frac{\partial w^1}{\partial x} + \frac{\partial w^2}{\partial y} \right) + \varepsilon^{22} \frac{\partial w^2}{\partial x} \right] \zeta ~d\Omega 
        - 
        \frac{1}{2}\int_{\Omega(t)} \left( \frac{\partial w^1}{\partial y} + \frac{\partial w^2}{\partial x} \right) \zeta ~d\Omega.
    \end{align*}
    \item If $\boldsymbol{\zeta} = \begin{pmatrix} 0 & 0 \\ 0 & \zeta \end{pmatrix}$, then 
    \begin{align*}
        &  \int_{\Omega(t)} \left[ \boldsymbol{\varepsilon}\,\text{skw}(\nabla {\bf w}) - \text{skw}(\nabla {\bf w})\,\boldsymbol{\varepsilon} 
        + 
        \left( \text{tr}(\boldsymbol{\varepsilon})-1 \right)\,\text{sym}(\nabla {\bf w})
        - 
        (\nabla \cdot {\bf w}) \boldsymbol{{\varepsilon}} \right] : \boldsymbol{\zeta} ~ d\Omega = \\
        & \qquad
        = \int_{\Omega(t)} \left[ \varepsilon^{11} \dfrac{\partial w^2}{\partial y} - \varepsilon^{12} \left( \frac{\partial w^2}{\partial x} - \frac{\partial w^1}{\partial y} \right) - \varepsilon^{22} \frac{\partial w^1}{\partial x} \right] \zeta ~d\Omega
        -
        \int_{\Omega(t)} \frac{\partial w^2}{\partial y} \zeta ~d\Omega.
    \end{align*}
\end{itemize}
\end{rem}

Let $\mathcal{T}_h(0)$ be a regular triangulation of the initial domain $\Omega(0)$. For spatial discretization, we consider a conforming finite element approximation on a time-dependent triangulation $\mathcal{T}_h(t)$. We introduce the discrete spaces
\[
{\mathcal{W}}_h \subset (H^1(\Omega))^d, \quad 
{\mathcal{Q}}_h \subset H^1(\Omega), \quad 
{\mathcal{E}}_h \subset (L^2(\Omega))^{d \times d},
\]
where ${\mathcal{W}}_h$ and ${\mathcal{Q}}_h$ consist of continuous piecewise linear (P1) finite elements. Further, ${\mathcal{E}}_h$ may consist of piecewise constant basis functions, however, currently we use piecewise linear (P1) elements for ${\mathcal{E}}_h$ in our code.
For the Galerkin finite-element formulation with time-dependent linear basis functions we write 
\begin{align*}
	\mathbf{w}_h(\mathbf{x},t) = \sum_{k=1}^{n_w} \mathbf{w}_k(t) \boldsymbol{\phi}_k(\mathbf{x},{\bf a}(t)), 
    \quad    \boldsymbol{\varepsilon}_h(\mathbf{x},t) = \sum_{k=1}^{n_\varepsilon} \boldsymbol{\varepsilon}_k(t) \boldsymbol{\zeta}_k(\mathbf{x},{\bf a}(t)),
    \quad 
	\boldsymbol{p}_h(\mathbf{x},t) = \sum_{k=1}^{n_p} p_k(t) \psi_k(\mathbf{x},{\bf a}(t)),
\end{align*}
where ${\bf a}(t)$ denotes the element coordinates (determined by displaced vertices) in time, $n_w$, $n_\varepsilon$ and $n_p$ denote, respectively, the number of degrees of freedom for the discrete velocity, strain and pressure fields. 

Substitution of the Galerkin approximations in the weak form results in the following semi-discrete algebraic system: 
\begin{subequations}  \label{mvpe_semidisc}
\begin{align}
	& \rho \frac{d}{dt} (M_w {\underline w}) + S_{vis} {\underline w} + S_{el} {\underline \varepsilon} - D^T {\underline p} = {\underline b}_w, \\[1ex]
	%& \frac{d}{dt} (M_{\varepsilon} {\underline \varepsilon})  + N({\underline w},{\underline \varepsilon}) = - \underline{b}_G, \\[1ex]
    & 
    \frac{d}{dt} (M_{\varepsilon} {\underline \varepsilon}) - B \underline{w} + N({\underline w},{\underline \varepsilon}) + \alpha M_{\varepsilon} \underline{\varepsilon} = \underline{0}, \\[1ex]
	& D {\underline w} + \kappa L {\underline p} = %{\underline 0}.
    \underline{b}_p.
\end{align}
\end{subequations}
Here, vectors $\underline{b}_w$ and $\underline{b}_p$ account for the right hand side terms in Eqs. \eqref{wf_w} and \eqref{wf_p} respectively, $M_w$ and $M_\varepsilon$ the (block) mass matrices, $D$ denotes the divergence matrix, $S_{el}$ and $S_{vis}$ the elasticity and viscous stiffness matrices, $L$ the Laplace matrix and matrix $B$ and vector $N(\underline{w},\underline{\varepsilon})$ arise from the discretization of the central term in Eq. \eqref{wf_varepsilon}, see Remark \ref{rem:wf_eps}.

Let $\Delta t$ be the time step and denote by $\tau$ the current time level (i.e., $t^\tau = \tau \Delta t$). After application of the backward Euler method to the semi-discrete problem in Eqs. (\ref{mvpe_semidisc}), one obtains 
\begin{subequations}\label{discrmorphporo1}
\begin{align}
	& \Big( \! \rho M_w^{\tau} + \Delta t S_{vis}^{\tau} \! \Big)  {\underline w}^{\tau} + \Delta t S_{el}^{\tau} {\underline \varepsilon}^{\tau} - \Delta t (D^{\tau})^T {\underline p}^{\tau} = \rho M_w^{\tau-1} {\underline w}^{\tau-1} +  \Delta t {\underline b}_w^{\tau},  \label{discrmorphporo1_a}\\[1ex]
	& (1+\Delta t) M_\varepsilon^{\tau} \underline{\varepsilon}^{\tau} - \Delta t B^{\tau} \underline{w}^{\tau} + \Delta t \ N({\underline w}^{\tau},{\underline \varepsilon}^{\tau})  = M_\varepsilon^{\tau-1} {\underline \varepsilon}^{\tau-1}, \label{discrmorphporo1_b}\\[1ex] %G^{\tau}, \\[1ex]
	& D^{\tau} {\underline w}^{\tau} + \kappa L^{\tau} {\underline p}^{\tau} = \underline{b}^\tau_p,
    \label{discrmorphporo1_p}
\end{align}
\end{subequations}
which contains time-dependent matrices and vectors. Note that Eq. (\ref{discrmorphporo1_b}) yields a nonlinear equation for the strain $\underline{\varepsilon}$. We adopt a fixed point iteration, on each time level, to solve the resulting nonlinear problem.
As we shall see in the numerical results, non-physical oscillations may appear in the pressure field for certain cases and later we analyze the conditions that lead to these oscillations in the one-dimensional case. 

\subsection{Monotonicity requirements} \label{sec:monotonicity}
We consider the one-dimensional problem and, for notational convenience, we drop the $\tau$ superscripts on the matrices and vectors at time $t^\tau$. After algebraic manipulation of Eqs. (\ref{discrmorphporo1_a}) and (\ref{discrmorphporo1_p}), we obtain that the pressure $\underline{p}$ must satisfy the following equation
\begin{equation}
A \underline{p} = \Delta t D P^{-1} S_{el}\underline{\varepsilon} + \tilde{\underline{b}}_p,
\end{equation}
where 
\[
    A = \kappa L + \Delta t D P^{-1}D^T, \qquad P = \rho M + \Delta t S_{vis} = \rho M + \Delta t \mu_{vis} L,
\]
with $\mu_{vis} = \mu_1 + \mu_2$, and 
\[
    \tilde{b}_p = -\rho D P^{-1}M^{\tau-1} \underline{w}^{\tau-1} - \Delta t D P^{-1} \underline{b}_w.
\]
If $A$ is an $M$-matrix then the pressure field satisfies a discrete maximum principle, which prohibits the occurrence of non-physical oscillations. Here, $P^{-1}$ is approximated by using fundamental solutions \cite{Vermolen2022}, which results in the following approximation 
\begin{subequations} \label{approx}
\begin{align} 
	& (P^{-1})_{ij} \approx (\tilde{P}^{-1})_{ij} = \frac{\nu}{2 \rho \cosh(\nu)} \bar{p}_{ij} \\
	& \bar{p}_{ij} = 
	\begin{cases}
		\sinh(\nu(1-x_i+x_j)) + \sinh(\nu(1-x_i-x_j)) , & \text{if $j \le i$,} \\[1ex]
		\sinh(\nu(1+x_i-x_j)) + \sinh(\nu(1-x_i-x_j)) , & \text{if $j > i$,}
	\end{cases}
\end{align}
\end{subequations}
where $x_i$ represents the nodal point $i$ and
\[
	\nu := \sqrt{\dfrac{\rho}{\Delta t \mu_{vis}}}.
\]
In \cite{Vermolen2022}, it has been proved that $| (P^{-1})_{ij} - (\tilde{P}^{-1})_{ij} | = \mathcal{O}(h^{3/2})$, where $h$ denotes the size of the finite-element mesh. 

We note that the analysis can be done for a non-uniform, however, for the sake of illustration, we give the analysis for a uniform mesh. 
 
After some tedious algebraic manipulations, see Section $4.4$ in \cite{denBakker2020}, it follows that
\begin{equation*}
(D P^{-1} D^T {\underline p})_i = \frac{h}{4  \Delta t \mu_{vis} } (p_{i-1} + 2 p_i + p_{i+1}) + \mathcal{O}(h^{7/2}),
\end{equation*}
since the coefficients of $D$ are of order $h$. Neglecting higher order terms in the previous construction, we obtain

\begin{align*}
	(A {\underline p})_i & = \left( (\kappa L + \Delta t DP^{-1}D^T) p \right)_i \\
    & \approx \left(  \frac{h}{4\mu_{vis}} - \frac{\kappa}{h} \right) p_{i-1} + 2\left( \frac{h}{4\mu_{vis}} + \frac{\kappa}{h} \right) p_i + \left( \frac{h}{4\mu_{vis}} - \frac{\kappa}{h} \right) p_{i+1}.
\end{align*}

Recall that for a matrix to be an $M$--matrix (for monotonicity of the finite-element solution), the off-diagonals need to be non--positive. %, hence we need

\begin{thm}\label{Theorem1}
The matrix $A$ based on approximation (\ref{approx}) for $P^{-1}$ is an $M$-matrix if
\begin{equation} \label{eq:h_constrain} 
h \leq 2\sqrt{\mu_{vis}\kappa}.
\end{equation}
We remark that this condition is sufficient to remove spurious oscillations for the approximate system. If $\mu_{vis}   \rightarrow 0$, then the above result becomes sharper. The assertion displays an indication that (large) violation of this condition leads to spurious oscillations. We will observe this behavior in the numerical results section.
\end{thm}

Note that this constraint affects only the mesh size and not the time step. 
If the viscosity and/or permeability are low, then the mesh resolution has to be very large. For these cases, stabilization may be attractive. On the other hand, note also that for $\kappa \longrightarrow \infty$ the above criterion will always be satisfied.

\begin{rem}
{As $ \kappa \rightarrow \infty$, then the morphoelastic model is recovered, which always gives oscillation-free finite element approximations.}
\end{rem}

\subsection{Stabilized finite-element method}

We follow the stabilization technique used by Aguilar et al. \cite{Aguilar2008} and Rodrigo \cite{RodrigoETAL2016}, which is based on the perturbation of the pressure equation. Note that the result that we are deriving in this section holds for the one-dimensional case. In our case, we introduce this perturbation for the discrete system of equations. Thus, we replace Eq. (\ref{discrmorphporo1_p}) by 
\begin{equation} \label{discrmorphporo}  
  D^{\tau} {\underline w}^{\tau} + (\kappa + \beta) L^{\tau} {\underline p}^{\tau} = \beta L^{\tau-1} \underline{p}^{\tau-1} + \underline{b}^\tau_p, 
\end{equation}
where $\beta>0$ is the stabilization parameter. For the sake of presentation, we drop the superscript $\tau$. We look for an appropriate choice of $\beta$ such that 
\[
    B = (\kappa + \beta)L +  \Delta t D P^{-1}D^T
\]
is an $M$--matrix. Using the approximation to $P^{-1}$ introduced in Section \ref{sec:monotonicity} we have, for the one-dimensional problem, that
\begin{align*}
	(B {\underline p})_i & = \left( \Big((\kappa+\beta) L + \Delta t DP^{-1}D^T\Big) p \right)_i \\
    & \approx \left(  \frac{h}{4\mu_{vis}} - \frac{\kappa+\beta}{h} \right) p_{i-1} + 2\left( \frac{h}{4\mu_{vis}} + \frac{\kappa+\beta}{h} \right) p_i + \left( \frac{h}{4\mu_{vis}} - \frac{\kappa+\beta}{h} \right) p_{i+1}.
\end{align*}

The approximate matrix is an $M$--matrix if and only if
\[
    \frac{h}{4\mu_{vis}}-\frac{\kappa + \beta}{h} \leq 0 
    \qquad \iff \qquad 
    \frac{h^2}{4\mu_{vis}} - \kappa \leq \beta.
\]
\begin{thm}\label{Theorem4}
{The matrix $B$ based on approximation (\ref{approx}) for $P^{-1}$ is an $M$-matrix if
\begin{equation}
\beta \ge \beta^* = \frac{h^2 }{4 \mu_{vis} }-\kappa,
\end{equation}

A sufficient condition for this is
$\beta \ge \dfrac{ h^2 }{4 \mu_{vis} }$.
}

\settowidth{\leftmargin}{{\Large$\square$}}\advance\leftmargin\labelsep
\itemsep8pt\relax
\renewcommand\labelitemi{{\lower1.5pt\hbox{\Large$\square$}}}
\end{thm} 

The above condition is based on small $h$, and in the limit, warrants monotonicity of the finite element solution regardless of the permeability parameter $\kappa$ and the time step $\Delta t$. 
We expect a similar condition for higher dimensionality in the future.

In our numerical simulations we solve Eqs. (\ref{discrmorphporo1_a}), (\ref{discrmorphporo1_b}), and (\ref{discrmorphporo}) to compute the solid matrix velocity $\underline{w}^\tau$, the components of the strain tensor $\underline{\varepsilon}^\tau$, and the pressure $\underline{p}^\tau$ at time level $\tau$. Note that Eq. (\ref{discrmorphporo1_b}) contains a nonlinear coupling between the velocity and the strain. We employ a fixed-point iteration method combined with stabilization at each time step. The used algorithm is described in Algorithm~\ref{algo}.
 
\begin{algorithm}
\caption{Morpho-Visco-Poroelastic FEM Solver (Fixed-point Iteration).}
\label{algo}
\begin{algorithmic}[1]
\REQUIRE
    Initial triangulation $\mathcal{T}_h(0)$ of $\Omega(0)$,
    initial conditions for $\underline{w}$, $\underline{\varepsilon}$, value of
    parameters $\rho, \mu_{vis}, \kappa, \alpha, \beta$,
    time step $\Delta t$, final time $T$,
    source terms $\underline{b}_w^\tau$, $\underline{b}_p^\tau$ for each $\tau$.
\ENSURE
    Discrete solutions $\underline{w}^\tau$, $\underline{\varepsilon}^\tau$, $\underline{p}^\tau$ for $\tau = 1, 2, \dots, T/\Delta t$.

\FOR{$\tau = 1$ to $T/\Delta t$}
    \STATE Update mesh coordinates $\mathbf{a}(t^\tau)$ from displacement field.
    \STATE Assemble time-dependent matrices on the updated configuration:
        $M_w^\tau$, $S_{vis}^\tau$, $S_{el}^\tau$, $D^\tau$, $M_\varepsilon^\tau$, $L^\tau$, $B^\tau$.
    
    \STATE Initialize fixed-point iteration: $k = 0$, $\underline{\varepsilon}^{\tau,0} = \underline{\varepsilon}^{\tau-1}$.
    
    \REPEAT
        \STATE \textbf{Step 1:} Solve for velocity and pressure from the linear system
        \[
        \begin{cases}
        \bigl( \rho M_w^\tau + \Delta t S_{vis}^\tau \bigr) \underline{w}^{\tau,k+1}
        + \Delta t S_{el}^\tau \underline{\varepsilon}^{\tau,k}
        - \Delta t (D^\tau)^T \underline{p}^{\tau,k+1}
        = \rho M_w^{\tau-1} \underline{w}^{\tau-1} + \Delta t \underline{b}_w^\tau, \\[2ex]
        D^\tau \underline{w}^{\tau,k+1} + (\kappa + \beta) L^\tau \underline{p}^{\tau,k+1}
        = \beta L^{\tau-1} \underline{p}^{\tau-1} + \underline{b}_p^\tau .
        \end{cases}
        \]
        
        \STATE \textbf{Step 2:} Update strain by solving the linearized equation
        \[
        (1+\Delta t) M_\varepsilon^\tau \underline{\varepsilon}^{\tau,k+1}
        - \Delta t B^\tau \underline{w}^{\tau,k+1}
        + \Delta t \, N(\underline{w}^{\tau,k+1}, \underline{\varepsilon}^{\tau,k})
        = M_\varepsilon^{\tau-1} \underline{\varepsilon}^{\tau-1}.
        \]
        
        \STATE \textbf{Step 3:} Compute residual $\delta = \|\underline{\varepsilon}^{\tau,k+1} - \underline{\varepsilon}^{\tau,k}\|$.
        \STATE Set $k \leftarrow k+1$.
    \UNTIL{$\delta < \text{tol}$}
    
    \STATE Set $\underline{w}^\tau = \underline{w}^{\tau,k}$, $\underline{\varepsilon}^\tau = \underline{\varepsilon}^{\tau,k}$, $\underline{p}^\tau = \underline{p}^{\tau,k}$.
\ENDFOR

\RETURN $\{\underline{w}^\tau, \underline{\varepsilon}^\tau, \underline{p}^\tau\}_{\tau=1}^{T/\Delta t}$.
\end{algorithmic}
\end{algorithm}

\section{Numerical results} \label{sec:5}
In this section we consider a two dimensional test problem. We test the numerical solution on the morpho-visco-poroelastic model, Eqs. \eqref{morphoporoelas} in the unit square,  $\Omega(0) = (0,1)^2$, subject to homogeneous boundary conditions 
\begin{subequations}
  \begin{empheq}[left=\empheqlbrace]{align}
   & {\bf w} = \mathbf{0}, \quad \kappa \nabla p \cdot \boldsymbol{n} = \mathbf{0} , & \text{on $\Gamma_1(t) \times (0,T]$,} \\[1ex]
   & \boldsymbol{\sigma}\cdot\boldsymbol{n} = \mathbf{0}, \quad p = 0, & \text{on $\Gamma_2(t) \times (0,T]$,}
  \end{empheq}
\end{subequations}
where $\Gamma_1 = \{{(x, y) : x = 0, y \in [0, 1]\}}$ and $\Gamma_2 = \partial \Omega \setminus \Gamma_1$, and initial conditions are given in Eqs. \eqref{mpve_ic1}-\eqref{mpve_icw}. We consider the problem with constant mechanical parameters
as given in Table \ref{param}.
\begin{table}[htbp]
\centering
\caption{Parameter values for the numerical simulations in Section \ref{sec:5}.} \label{param}
\begin{tabular}{cccccccc}
\hline
$\mu$ & $\lambda$ & $\mu_{vis}$ & $\rho$ &  $\alpha$ \\
\hline
0.5 & 1 & 1 & 1 & 1 \\
\hline
\end{tabular}
\end{table}
In addition, the material is subjected to a time-dependent body force ${\bf f}_{\bf u}$, see Eq. \eqref{mvpe_momentum}, given by 
\begin{equation}
{\bf f}_{\bf u}(t,{\bf x}) = 
\begin{pmatrix}
0 \\
e^{-t} \sin(2\pi t)
\end{pmatrix}.
\end{equation}

\subsection{Pressure profiles and stabilization}
\begin{figure}[h!]
    \centering
    % First row
    \begin{subfigure}[b]{0.495\textwidth}
        \centering
        \includegraphics[width=\textwidth]{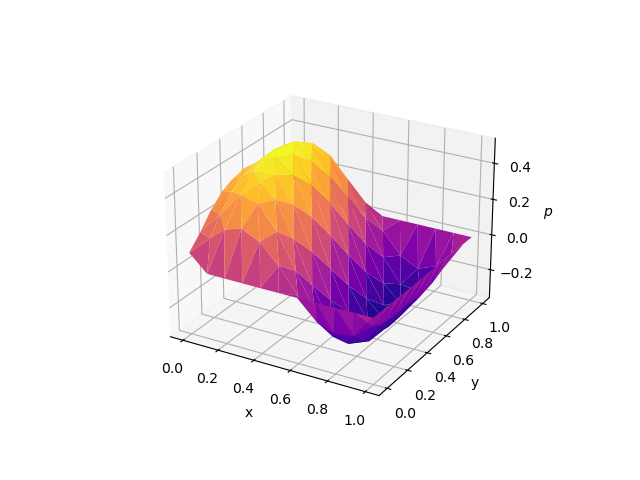}
        \caption{$\kappa=10^{-2}$ and $h=0.14$.}
        \label{fig:non_stabilized_1}
    \end{subfigure}
    \hfill
    \begin{subfigure}[b]{0.495\textwidth}
        \centering
        \includegraphics[width=\textwidth]{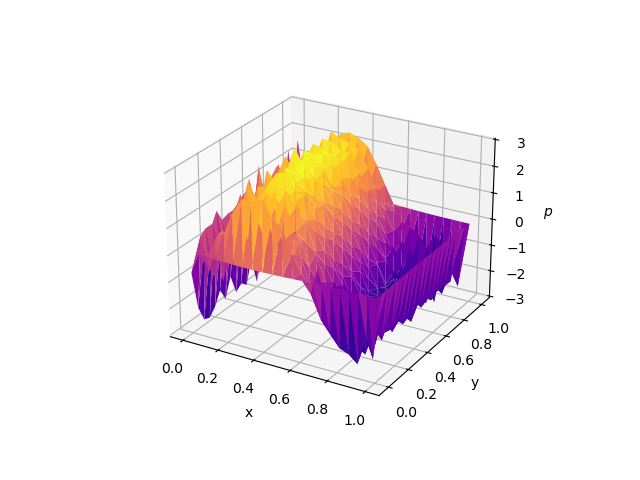}
        \caption{$\kappa=10^{-6}$ and $h=0.07$.}
        \label{fig:non_stabilized_2}
    \end{subfigure}
    \caption{Pressure fields for different pairs of $\kappa$ and $h$.}
    \label{fig:pressurefield_non_stabilized}
\end{figure}

Fig. \ref{fig:pressurefield_non_stabilized} shows the pressure field after one time step, $\Delta t = 0.1$, for different pairs of permeability $\kappa$ and mesh size $h$ (the diameter of the triangular elements). We can appreciate that the pressure field for $\kappa=10^{-2}$ and $h=0.14$ remains smooth (Fig. \ref{fig:non_stabilized_1}). This is to be expected as the chosen permeability and mesh size fulfill Eq. \eqref{eq:h_constrain}, $h \leq 2\sqrt{\mu_{vis}\kappa}$. However, the pair $\kappa=10^{-6}$ and $h=0.07$ does not fulfill Eq. \eqref{eq:h_constrain} and spurious oscillations are visible in the pressure field (Fig. \ref{fig:non_stabilized_2}).

\begin{figure}[h!]
    \centering
    % First row
    \begin{subfigure}[b]{0.32\textwidth}
        \centering
        \includegraphics[width=\textwidth]{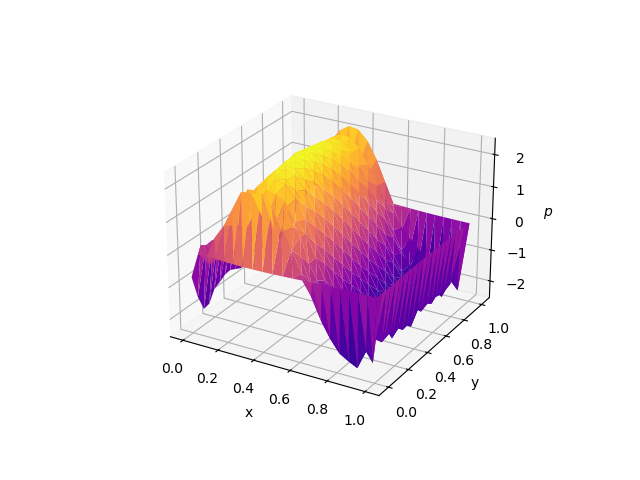}
        \caption{$\beta = 10^{-5}$.}
        \label{fig:stabilized_1}
    \end{subfigure}
    \hfill
    \begin{subfigure}[b]{0.32\textwidth}
        \centering
        \includegraphics[width=\textwidth]{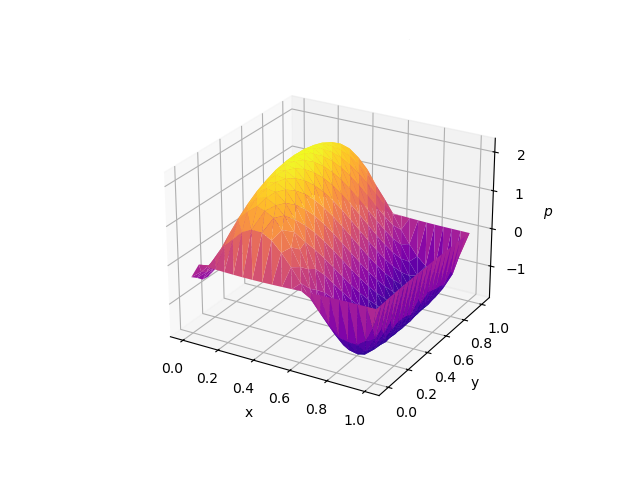}
        \caption{$\beta = 1.25 \times 10^{-3}$.}
        \label{fig:stabilized_2}
    \end{subfigure}
    \hfill
    \begin{subfigure}[b]{0.32\textwidth}
        \centering
        \includegraphics[width=\textwidth]{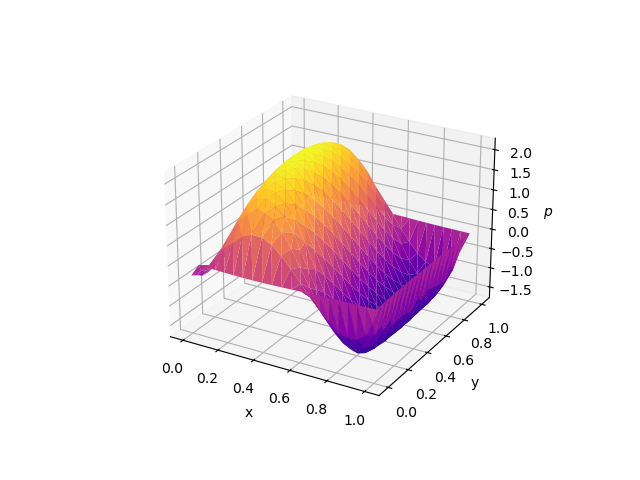}
        \caption{$\beta = 10^{-3}$.}
        \label{fig:stabilized_3}
    \end{subfigure}
    \caption{Pressure fields for $\kappa = 10^{-6}$ and $h=0.07$ and different values of the stabilization parameter $\beta$.}
    \label{fig:pressurefield_stabilized}
\end{figure}

In order to obtain a monotonic pressure field, we could increase the grid resolution to meet Eq. \eqref{eq:h_constrain} (which yields $h\leq 0.0028$). We could also apply the stabilization technique proposed in the previous section. Theorem \ref{Theorem4} establishes that the approximated system will stabilize for $\beta \geq 1.25\times 10^{-3}$. Fig. \ref{fig:pressurefield_stabilized} shows the pressure profiles for different values of the stabilization parameter $\beta$. We observe the correspondence of the numerical results with the theoretical estimate. 

\subsection{Total Variation ($\mathrm{TV}$) for pressure field}
\begin{figure}[h!]
\centering
\includegraphics[width=10cm]{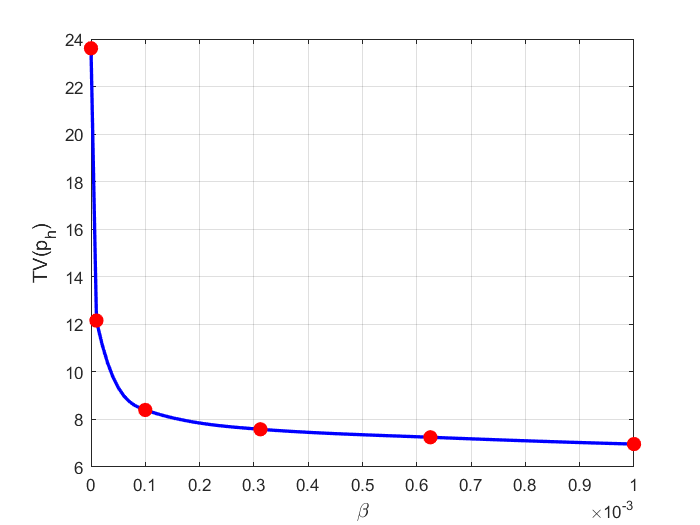}
\caption{\textit{$\mathrm{TV}$ for pressure fields for $\kappa = 10^{-6}$ and $h=0.07$ versus different values of the stabilization parameter $\beta$ for one time step $\Delta t=0.1$}.}
\label{TV_beta}
\end{figure}

In order to have an objective measure for the amount of variation of the numerical solution, we consider the total variation over the domain of computation. In our case $\Omega$ is a rectangular domain. Let the numerical solution of the pressure be given by $p_h$ and suppose that the horizontal and vertical sides of the rectangular triangular elements have lengths $\Delta x$ and $\Delta y$, respectively. The discrete function $p_h$ is defined on this grid with $p_{i,j}$ being the value associated with a gridpoint ${{\bf x}_{i,j}}$. We use the definition for the Total Variation (TV) of Goodman and LeVeque \cite{Goodman1985}, which is given by:
\begin{equation}
   \text{TV}_{\tau}(p_h) = \sum_{(i,j)\in {\tau}} \left( \Delta y |p_{i+1,j} - p_{i,j}| + \Delta x |p_{i,j+1} - p_{i,j}| \right).
   \label{TV_def}
\end{equation}
This definition is commonly referred to as anisotropic TV in the literature, where $\tau$ is the set of all mesh points.
This formula approximates the given continuous total variation functional for pressure field given by
\begin{equation}
\mathrm{TV}(p) = \int_\Omega |p_x| + |p_y| \, dxdy.
\end{equation}
We subsequently experimentally demonstrate that our stabilization technique is Total Variation Diminishing (TVD) with increasing $\beta$-values i.e. $$0 \le \mathrm{TV}_{\tau}(p_{\beta_{k+1}}) \leq \mathrm{TV}_{\tau}(p_{\beta_{k}}), ~ \forall ~ \beta_{k+1} \ge \beta_k.$$ 
We list the total variation of the numerical solutions for the pressure field using Eq. \eqref{TV_def} in Table~\ref{Table_TV_vs_beta}. We observe that the total variation of the solutions over the whole domain (i.e. unit square, $\Omega(0) = (0,1)^2$) monotonically decreases with increasing values of the stabilization parameter, which is to be expected (see also Fig. \ref{TV_beta}). Note that the total variation converges to a non-zero limit since the actual numerical solution is not constant over the domain of computation.
\begin{table}[h!]
\centering
\caption{\it Total variation ($\mathrm{TV}$) values for the pressure fields for different values of stabilization parameter $\beta$.} 
\label{Table_TV_vs_beta} 
\begin{tabular}{|c|c|c|}
\hline
Stabilization parameter ($\beta$) & Total Variation ($\mathrm{TV}(p_h)$) \\ 
\hline \hline
$0$ & 23.6096 \\ 
\hline
$10^{-5}$ & 12.1596 \\ 
\hline
$10^{-4}$ & 8.3987 \\ 
\hline
$3.12\times10^{-4}$ & 7.5860 \\ 
\hline
$6.25\times10^{-4}$ & 7.2501 \\
\hline
$10^{-3}$ & 6.9645 \\
\hline
\hline
\end{tabular}
\end{table}

\subsection{Convergence study}
Since no closed form solution is available, we analyze the convergence of the proposed discretization by Richardson extrapolation. We compute the $L^2(\Omega)$-- and $H^1(\Omega)$--norms of the velocity field ${\bf w}$ and the pressure $p$, using a fixed time step, at the final time $t_{end}$ of the simulation. The parameter values considered in this analysis are given in Table~ \ref{param}. For $\kappa = 10^{-2}$, oscillations are prevented provided that the element diameter satisfies $h\leq 0.2$, which holds in all the simulations that were shown in Table~ \ref{tab:beta=0} and hence stabilization is not required here.

The results show approximately second-order convergence in the \(L^2(\Omega)\)-norm for both \(\bf{w}\) and \(p\). In the \(H^1(\Omega)\)-norm, the observed rate is approximately first order for \(\mathbf{w}\), while the pressure attains approximately second-order convergence rate.
We obtained a similar convergence behavior when stabilization was applied, with $\beta = 10^{-3}$ (see Table ~\ref{tab:beta=e-3}), indicating that the proposed stabilization does not alter the convergence of the scheme and that the impact of the stabilization term is limited. 

Table~\ref{tab:time} gathers the errors obtained under successive refinements of the time step, when the finest mesh size is used and for final time \(t_{{end}}=0.1\). Again, no stabilization is required in this case (i.e. \(\beta=0\)).
The computed convergence orders are approximately between first and second order for both variables, demonstrating appropriate accuracy of the proposed discretization. Moreover, decreasing the time step improves the approximation quality, especially for the pressure variable $p$.

\begin{table}[htbp]
\centering
\caption{Spatial errors and convergence orders of $\bf w$ and $p$ with $\Delta t = 0.000625$, $t_{end}=0.0125$, and $\beta=0$.}
\label{tab:beta=0}
\begin{tabular}{m{1.5cm}|m{2.2cm}m{1cm}m{2.2cm}m{1.3cm}|m{1.5cm}m{1cm}m{2.1cm}m{1.3cm}}
\toprule
\multirow{2}{*}{$h/\sqrt{2}$} & \multicolumn{4}{c|}{$\bf w$} & \multicolumn{4}{c}{$p$} \\
\cmidrule(lr){2-5} \cmidrule(lr){6-9}
& $\lVert {\bf w} \rVert_{L^2}$ & order & $\lVert {\bf w} \rVert_{H^1}$ & order
& $\lVert p\rVert_{L^2}$ & order & $\lVert p\rVert_{H^1}$ & order \\
\midrule
$1/32$   & 2.104012e-2 &   1.9594    & 1.609170e-1 & 1.0033 & 8.889836  &   2.2359    & 1.951037e1 & 2.2724 \\
$1/16$   & 2.084597e-2 & 1.2972 & 1.522910e-1 & 0.5087& 8.844537  & 1.3296 & 1.946016e1  &  0.4116\\
$1/8$  & 2.009093e-2 &  &  1.349999e-1 &  & 8.631157  &  & 1.921759e1 &  \\
$1/4$  & 1.823547e-2&  & 1.103982e-1 &  & 8.094854  &  & 1.889494e1 &  \\
\bottomrule
\end{tabular}
\end{table}

\begin{table}[htbp]
\centering
\caption{Spatial errors and convergence orders of $\bf w$ and $p$ with $\Delta t = 0.000625$, $t_{end}=0.0125$, and $\beta=10^{-3}$.}
\label{tab:beta=e-3}
\begin{tabular}{m{1.5cm}|m{2.2cm}m{1cm}m{2.2cm}m{1.3cm}|m{1.5cm}m{1cm}m{2.1cm}m{1.3cm}}
\toprule
\multirow{2}{*}{$h/\sqrt{2}$} & \multicolumn{4}{c|}{$\bf w$} & \multicolumn{4}{c}{$p$} \\
\cmidrule(lr){2-5} \cmidrule(lr){6-9}
& $\lVert {\bf w} \rVert_{L^2}$ & order & $\lVert {\bf w} \rVert_{H^1}$ & order
& $\lVert p\rVert_{L^2}$ & order & $\lVert p\rVert_{H^1}$ & order \\
\midrule
$1/32$   & 2.013644e-2 &   1.9551    & 1.545243e-1 &   0.9868    & 8.906786  &   2.2308    & 1.954305e1 &   2.2615    \\
$1/16$   & 1.994895e-2 & 1.2952 & 1.460621e-1 & 0.4990 & 8.860838  & 1.3309 & 1.949115e1     &  0.4591\\
$1/8$  & 1.922199e-2 &  & 1.292921e-1  &  & 8.645168  &  & 1.924230e1  &  \\
$1/4$  & 1.743791e-2&  & 1.055924e-1 &  & 8.102622  &  & 1.890021e1  &  \\
\bottomrule
\end{tabular}
\end{table}

\begin{table}[htbp]
\centering
\caption{Spatial errors and convergence orders of $\bf w$ and $p$ with $h=0.03125\sqrt{2}$, $t_{end}=0.1$, and $\beta=0$.}
\label{tab:time}
\begin{tabular}{m{1.3cm}|m{2.2cm}m{1cm}m{2.2cm}m{1.3cm}|m{1.5cm}m{1cm}m{1.4cm}m{2.1cm}}
\toprule
\multirow{2}{*}{$\Delta t$} & \multicolumn{4}{c|}{$\bf w$} & \multicolumn{4}{c}{$p$} \\
\cmidrule(lr){2-5} \cmidrule(lr){6-9}
& $\lVert {\bf w} \rVert_{L^2}$ & order & $\lVert {\bf w} \rVert_{H^1}$ & order
& $\lVert p\rVert_{L^2}$ & order & $\lVert p\rVert_{H^1}$ & order \\
\midrule
$0.00625$   & 1.276693e-1 &   1.3624    & 6.104306e-1 & 1.5329 & 1.009416  &  1.7441 & 4.917025 &  1.7032   \\
$0.0125$   & 1.270581e-1 & 1.4946 & 6.091477e-1 & 1.7909 & 1.017158  & 2.0839 & 4.930888  & 2.0377 \\
$0.025$  & 1.254866e-1 &  1.3513&  6.054353e-1 & 1.6096 & 1.043093  & 1.9644 & 4.976030 & 1.9170 \\
$0.05$  & 1.210529e-1&  & 5.925889e-1 &  & 1.153042  &  &5.161384&  \\
$0.1$  & 1.097408e-1&  & 5.533863e-1 &  & 1.582125  &  &5.861358 &  \\
\bottomrule
\end{tabular}
\end{table}

\subsection{Qualitative aspects of the morpho-visco-poroelastic model}
To see how the model works and which aspects are brought into play by the various features in the model, we provide some additional numerical results in this subsection.

\subsubsection{Vector plot for velocity field}
\begin{figure}[h!]
\centering
    \includegraphics[width=0.6\linewidth]{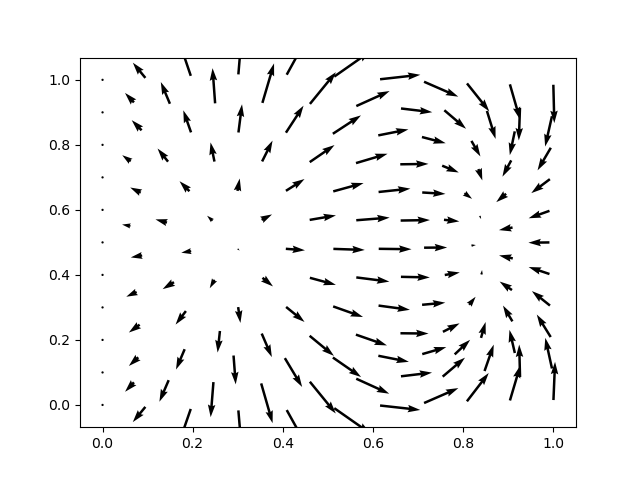}
\caption{Vector plot of the velocity field $\bf w$ of the morpho-visco-poroelastic model. Parameter values are taken from Table~\ref{param} with $\kappa = 10^{-4}$, hence, the stabilization term is not necessary here.}
\label{vector_plots}
\end{figure}

Fig. \ref{vector_plots} shows the solid velocity field in a unit square with respect to the parametric values given in Table~ \ref{param} with $\kappa=10^{-4}$. For the values in Table~ \ref{param}, the stabilization term is not necessary. Fig. \ref{vector_plots} is shown to illustrate what the solid displacement velocity looks like. It can be seen that through $\Gamma_1$ there is no flow, which is consistent with the boundary condition.  It can be seen that from the point $(0.3,0.5)$, the velocity shows a source with streamlines converging to a sink located at approximately $(0.9,0.5)$.

\subsubsection{Strain components}
\begin{figure}[htbp]
    \centering
        \includegraphics[width=1.02\linewidth]{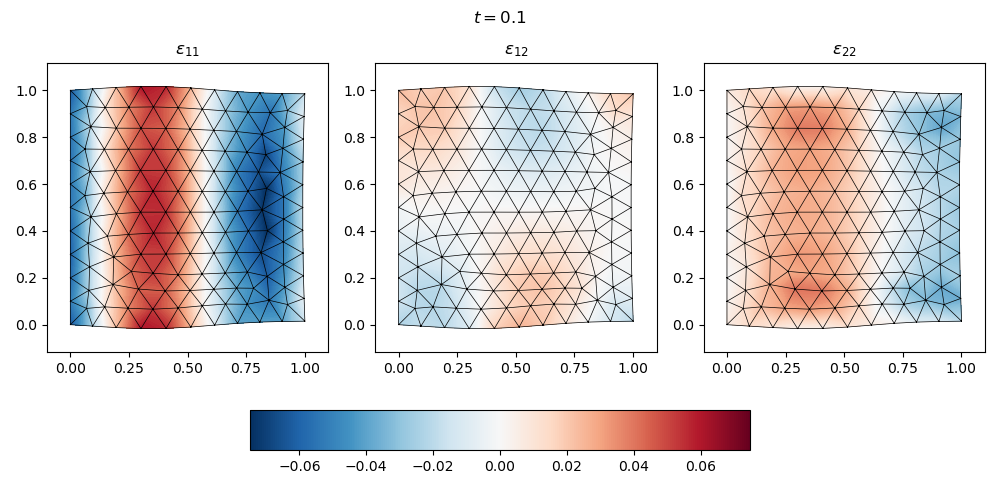}
    \caption{Strain components of $\boldsymbol{\varepsilon}$ in mopho-visco-poroelastic model with $h=0.14$ are shown. Parameter values are taken from Table~\ref{param}.}
    \label{figset1-non-stab-alpha}
\end{figure}

Fig. \ref{figset1-non-stab-alpha} shows the strain plots for the spatial distribution of the infinitesimal strain tensor components at $t=0.1$ over the deformed computational domain with respect to the parametric values given in Table~ \ref{param} with $\kappa=10^{-4}$. The plot mimics the impact of the combination of forces and microstructural changes. The component \(\varepsilon_{11}\) shows pronounced tensile and compressive regions along the horizontal direction, showing that the deformation is mainly governed by horizontal displacement gradients. The shear strain component \(\varepsilon_{12}\) exhibits comparatively smaller magnitudes, reflecting a mild angular distortion within the domain. Meanwhile, \(\varepsilon_{22}\) indicates moderate vertical deformation, with tensile behavior near the center and slight compression near the right boundary. In general, the strain magnitudes remain within a small range, consistent with the assumptions of infinitesimal strain theory.

\subsubsection{Comparison between different elastic models}
\begin{figure}[h!]
\includegraphics[width=0.5\linewidth]{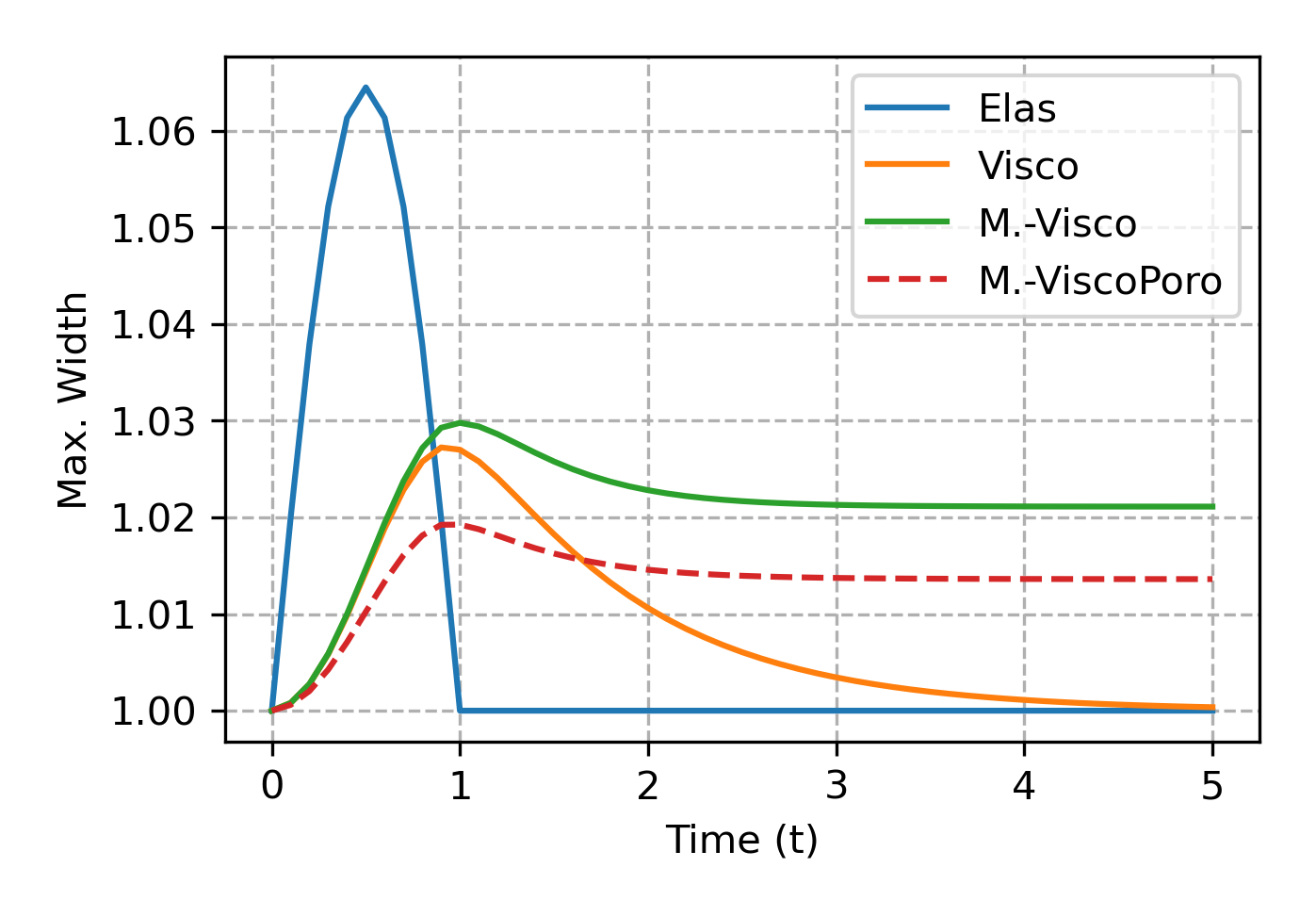}
\includegraphics[width=0.5\linewidth]{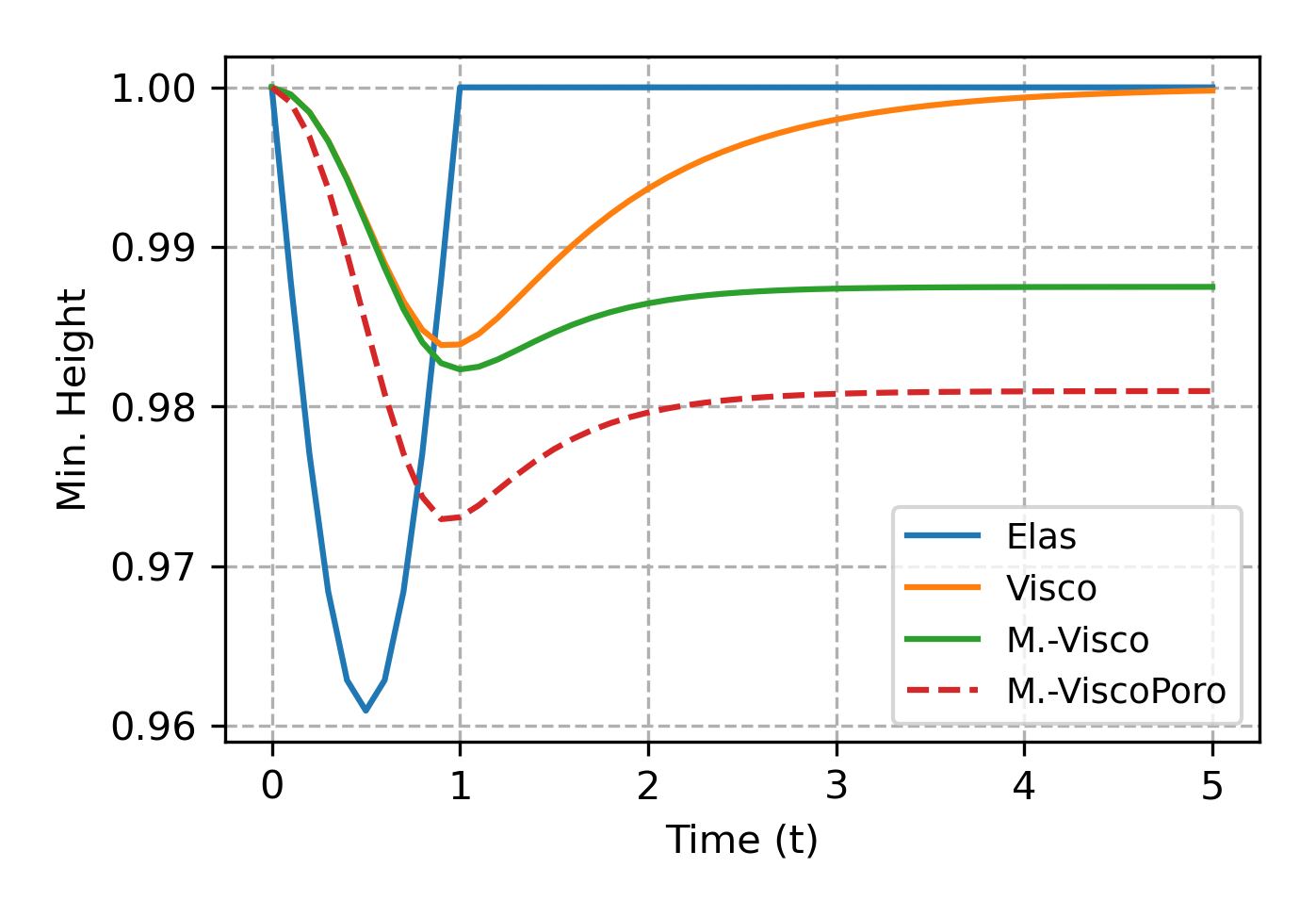}
\caption{A comparison of the time evolution of the maximal width and the minimal height for a horizontal stretching test for several models: elasticity, viscoelasticity, morpho-viscoelasticity, and morpho-visco-poroelasticity.}
\label{max_min}
\end{figure}

In this subsection, we give a qualitative comparison between four different elastic models. The values of the mechanical parameters are given in Table~ \ref{param} with $\kappa=10^{-4}$ and the time-dependent body force equal to:
\begin{equation*}
{\bf f_{u}}(t,{\bf x}) = 
\begin{pmatrix}
0.2 \sin(\pi t) \mathbbm{1}_{[0,1]} (t) \\ 
0
\end{pmatrix}.
\end{equation*}
Based on the above force, we expect the domain to stretch in the $x$-direction. To measure the stretching and compressing of the domain, we track the maximum domain width in some amount of transaxial compression as a consequence of axial stretching. Therefore, in order to monitor the stretching and compression of the domain, we track the maximum domain width in the horizontal direction and the minimal height in the vertical direction as a function of time. The simulations are run until $t_{end} = 50$ using a time-step of $\Delta t = 0.1$. The results are plotted in Fig. ~\ref{max_min}, where the evolution of the vertical $y$-coordinate is shown for different models: linear elasticity, visco-elasticity, morpho-viscoelasticity, and morpho-visco-poroelasticity. The linear elasticity model, where viscous effects and inertial effects have been disregarded, reacts instantaneously to the force. As soon as the force ceases, the configuration goes back to the initial state immediately. In this case, inertia has been neglected. The viscous case, which additionally incorporates the inertia term and a damping term due to viscous effects, bends back to its initial shape. In some simulations (not shown here), it oscillates slightly around the equilibrium during its return to the initial configuration as a result of the inertia term (which makes the equations behave like a (damped) hyperbolic wave equation). Furthermore, the viscous behavior, where an incompressible fluid is trapped in the pores, acts as a resistance to deformation, which makes the ultimate maximal deformation smaller than for the elastic case. For the morpho-viscoelastic case, it can be seen that the microstructural changes lead to the convergence to a revised equilibrium state that enhance larger displacement than for the visco-elastic case. Next to this model, we consider the morpho-visco-poroelastic formalism. The lowest peak deformation and the strongest damping effect are observed for this model. This behavior can be attributed to fluid-solid interaction and pressure dissipation within the porous medium.

\section{Discussion and Conclusions} \label{sec:6}
We have studied a morpho-visco-poroelastic model that combines the porous, elastic and viscous structure of tissues with microstructural changes that may lead to growth or shrinkage of the tissue. The model has been analyzed in terms of the symmetry of the strain tensor and linear stability. Linear stability was analyzed for both the continuous and semi-discrete systems. Stability conditions have been derived. In addition to stability, monotonicity of the numerical solution for the pressure can be warranted as long as the mesh size does not exceed a critical number that depends on the permeability of the porous medium. Furthermore, for larger mesh sizes, stabilization is required. This stabilization is based on the addition of a Laplacian operator to the pressure equation, and it has been found that the stabilization parameter should not be smaller than a threshold that incorporates the mesh size, viscosity, and permeability. The numerical simulations that we show confirm the stability bounds and the bounds for the stabilization that we derived theoretically.

Furthermore, our study was carried out for a one-dimensional case, and we want to extend this to higher dimensionality. For this purpose, we need expressions for the inverse of Laplace matrices and matrices that depend affinely on Laplace matrices for higher dimensions. We recently found some closed-form expressions for multi-dimensional matrix polynomials of Laplace matrices (with various boundary conditions) \cite{Asghar2026}, and we will study whether we can apply our recent findings in the current context.

The issue of well-posedness in terms of existence and uniqueness is an important matter that we want to consider in a separate manuscript. The proof will be inspired by the work of Both {\em et al.} \cite{Both2019}, who analyze gradient flow problems. In their interesting work, applications to thermo-visco-poroelasticity are studied.
 
\section*{Acknowledgments}
This work was supported by a fellowship awarded to SA by the Higher Education Commission (HEC) of Pakistan in the framework of project: 1(2)/HRD/OSS-III/BATCH-3/2022/HEC/527. The work of EJ is supported in part by the Spanish project PID2022-140108NB-I00 (MCIU/AEI/FEDER, UE), and by the DGA (Grupo de referencia APEDIF, ref. E24$\_$17R). The work of QP was supported by Research England under the Expanding Excellence in England (E3) funding stream, which was awarded to MARS: Mathematics for AI in Real-world Systems in the School of Mathematical Sciences at Lancaster University. 

\begin{appendices}
\section{}\label{Appendix}
\begin{table}[ht!]
\centering 
Nomenclature.
\begin{tabular}{p{1.75cm} p{8.75cm} p{4cm}}
\midrule
\textbf{Symbol} & \textbf{Description} & \textbf{Unit (SI)} \\
%\hline \hline
\midrule
\midrule
$\mathbf{w}$ & Displacement velocity & m$\cdot$s$^{-1}$ \\
$\boldsymbol{\varepsilon}$ & Strain tensor & dimensionless\\
$p$ & Pore fluid pressure & Pa (N$\cdot$m$^{-2}$)\\
$\mathbf{u}$ & Displacement vector %(post-processing) 
& m\\
$\rho$ & Tissue density & kg$\cdot$m$^{-3}$ \\
$\mu$ & Lamé parameter (shear modulus) & Pa\\
$\lambda$ & Lamé parameter (bulk modulus) & Pa\\
$E$ & Young's modulus ($E = 2\mu + \lambda$ in 1D) & Pa\\
$\mu_1, \mu_2$ & Viscosity coefficients (solid skeleton) & Pa$\cdot$s\\
$\mu_{vis}$ & Combined viscosity ($\mu_{vis} = \mu_1 + \mu_2$) & Pa$\cdot$s\\
$\kappa$ & Ratio of permeability and fluid viscosity & m$^3\cdot$s$\cdot${kg}$^{-1}$ (Pa$\cdot$s)$^{-1}$\\
$\alpha$ & Morphological growth/decay coefficient & s$^{-1}$\\
$\beta$ & Stabilization parameter & m$^3\cdot$s$\cdot${kg}$^{-1}$\\
$\mathbf{G}$ & Net growth tensor (rate of microstructural change) & s$^{-1}$\\
$\mathbf{f_u}$ & Applied body force vector & N$\cdot$m$^{-3}$ (kg$\cdot$m$^{-2}\cdot$s$^{-2}$)\\
$f_p$ & Fluid source/sink term & s$^{-1}$\\
$g_N$ & Fluid flux on boundary & m$\cdot$s$^{-1}$\\
$\mathbf{f}_b$ & External boundary force & Pa \\
${p}_0$ & Pressure constant & Pa \\
\hline \hline
\end{tabular}
\end{table}
\end{appendices}

\end{document}